\documentclass[11pt]{article}

\usepackage{amssymb}
\usepackage{amsmath}
\usepackage{amscd}
\usepackage[dvips]{graphicx}
\usepackage{mathrsfs}
\usepackage{wrapfig} 
\usepackage{picins}
\usepackage{picinpar}
\usepackage{psfig}
\usepackage{setspace}
\usepackage{float}

\numberwithin{equation}{section}

\newcommand{\tabincell}[2]{\begin{tabular}{@{}#1@{}}#2\end{tabular}}

\begin{document}

\title{On Spectral Properties of Finite Population\\ Processor Shared Queues}

\author{
Qiang Zhen\thanks{
Department of Mathematics and Statistics,
University of North Florida, 1 UNF Dr, Jacksonville, FL 32224, USA.
{\em Email:} q.zhen@unf.edu.}
\and
and
\and
Charles Knessl\thanks{
Department of Mathematics, Statistics, and Computer Science,
University of Illinois at Chicago, 851 South Morgan Street(M/C 249),
Chicago, IL 60607, USA.
{\em Email:} knessl@uic.edu.\
\newline\indent\indent{\bf Acknowledgement:} Q. Zhen was partly supported by a faculty scholarship grant from the University of North Florida. C. Knessl was partly supported by NSA grants H 98230-08-1-0102 and H 98230-11-1-0184. 
}}
\date{September 10, 2012}
\maketitle

\begin{abstract}
\noindent We consider sojourn or response times in processor-shared queues that have a finite population of potential users. Computing the response time of a tagged customer involves solving a finite system of linear ODEs. Writing the system in matrix form, we study the eigenvectors and eigenvalues in the limit as the size of the matrix becomes large. This corresponds to finite population models where the total population is $N\gg 1$. Using asymptotic methods we reduce the eigenvalue problem to that of a standard differential equation, such as the Hermite equation. The dominant eigenvalue leads to the tail of a customer's sojourn time distribution.   
\\

\noindent \textbf{Keywords:} Finite population, processor sharing, eigenvalue, eigenvector, asymptotics.
\end{abstract}

\section{Introduction}
The study of processor shared queues has received much attention over the past 45 or so years. The processor sharing (PS) discipline has the advantage over, say, first-in first-out (FIFO), in that shorter jobs tend to get through the system more rapidly. PS models were introduced during the 1960's by Kleinrock (see \cite{KLE_A}, \cite{KLE_T}). In recent years there has been renewed attention paid to such models, due to their applicability to the flow-level performance of bandwidth-sharing protocols in packet-switched communication networks (see \cite{HEY}-\cite{NAB}). 

Perhaps the simplest example of such a model is the $M/M/1$-PS queue. Here customers arrive according to a Poisson process with rate parameter $\lambda$, the server works at rate $\mu$, there is no queue, and if there are $\mathcal{N}(t)(>0)$ customers in the system each gets an equal fraction $(=1/\mathcal{N}(t))$ of the server. PS and FIFO models differ significantly if we consider the ``sojourn time''. This is defined as the time it takes for a given customer, called a ``tagged customer'', to get through the system (after having obtained the required amount of service). The sojourn time is a random variable that we denote by $\mathcal{V}$. For the simplest $M/M/1$ model, the distribution of $\mathcal{V}$ depends on the total service time $\mathcal{X}$ that the customer requests and also on the number of other customers present when the tagged customer enters the system.

One natural variant of the $M/M/1$-PS model is the finite population model, which puts an upper bound on the number of customers that can be served by the processor. The model assumes that there are a total of $N$ customers, and each customer will enter service in the next $\Delta t$ time units with probability $\lambda_0\Delta t+o(\Delta t)$. At any time there are $\mathcal{N}(t)$ customers being served and the remaining $N-\mathcal{N}(t)$ customers are in the general population. Hence the total arrival rate is $\lambda_0[N-\mathcal{N}(t)]$ and we may view the model as a PS queue with a state-dependent arrival rate that decreases linearly to zero. Once a customer finishes service that customer re-enters the general population. The service times are exponentially distributed with mean $1/\mu$ and we define the traffic intensity $\rho$ by $\rho=\lambda_0N/\mu$. This model may describe, for example, a network of $N$ terminals in series with a processor-shared CPU. This may be viewed as a closed two node queueing network.

The finite population model does not seem amenable to an exact solution. However, various asymptotic studies have been done in the limit $N\to\infty$, so that the total population, or the number of terminals, is large. If $N$ is large it is reasonable to assume either that $\lambda_0$, the arrival rate of each individual customer, is small, of the order $O(N^{-1})$, or that the service rate $\mu$ is large, of the order $O(N)$. Then $\rho=\lambda_0N/\mu$ will remain $O(1)$ as $N\to\infty$. Previous studies of the finite population model were carried out by Morrison and Mitra (see \cite{MIT}-\cite{MOR_C88}), in each case for $N\to\infty$. For example, the moments of the sojourn time $\mathcal{V}$ conditioned on the service time $\mathcal{X}$ are obtained in \cite{MOR_A}, where it was found that the asymptotics are very different according as $\rho<1$ (called ``normal usage''), $\rho-1=O(N^{-1/2})$ (called ``heavy usage''), or $\rho>1$ (called ``very heavy usage''). In \cite{MOR_M} the unconditional sojourn time distribution is investigated for $N\to\infty$ and the three cases of $\rho$, in \cite{MOR_C87} the author obtains asymptotic results for the conditional sojourn time distribution, conditioned on the service time $\mathcal{X}$, in the very heavy usage case $\rho>1$, and in \cite{MOR_C88} the results of \cite{MOR_C87} are generalized to multiple customer classes (here the population $N$ is divided into several classes, with each class having different arrival and service times). In \cite{MIT} the authors analyze the multiple class model and obtain the unconditional sojourn time moments for $N\to\infty$ in the normal usage case, while in \cite{MOR_H} heavy usage results are obtained.

In this paper we study the spectral structure of the finite population model as $N\to\infty$. 
We denote the sojourn time by $\mathcal{V}=\mathcal{V}(N)$ and its conditional density we call $p_n(t)$ with 
\begin{equation*}\label{s2_fs_pnt}
p_n(t)dt=\Pr\Big[\mathcal{V}(N)\in(t,t+dt)\Big|\mathcal{N}(0^-)=n\Big].
\end{equation*}
Here $\mathcal{N}(0^-)$ denotes the number of other customers present in the system immediately before the tagged customer arrives, and thus $0\le \mathcal{N}(0^-)\le N-1$. 
Then we define the column vector $\mathbf{p}(t)=(p_0(t),p_1(t),...,p_{N-1}(t))^T$. $\mathbf{p}(t)$ satisfies a system of ODEs in the form $\mathbf{p}'(t)=\mathbf{A}\mathbf{p}(t)$ where $\mathbf{A}$ is an $N\times N$ tridiagonal matrix, whose entries depend on $\rho=\lambda_0N/\mu$ and $N$. The eigenvalues of $\mathbf{A}$ are all negative and we denote them by $-\nu_j$ $(j=0,1,..., N-1)$ with the corresponding eigenvectors being $\phi_j(n)=\phi_j(n;N,\rho)$. We shall study this eigenvalue problem for $N\to\infty$ and three cases of $\rho$: $\rho<1$, $\rho>1$ and $\rho-1=O(N^{-1/3})$. In each case we obtain expansions of the $\nu_j$ and then the $\phi_j(n)$, for various ranges of $n$. Often the eigenvectors can be expressed in terms of Hermite polynomials for $N\to\infty$. Since $\mathbf{A}$ is a finite matrix the spectrum is purely discrete, but as the size of the matrix becomes large we sometimes see the eigenvalues coalescing about a certain value. Ordering the eigenvalues as $\nu_0<\nu_1<\nu_2<...$, the tail behavior of $p_n(t)$ and $p(t)$ for $t\to\infty$ is determined by the smallest eigenvalue $\nu_0$, where $p(t)$ is the unconditional sojourn time density with 
\begin{equation}\label{s1_pt}
p(t)=\sum_{n=0}^{N-1}p_n(t)\mathrm{Pr}\big[\mathcal{N}(0^-)=n\big].
\end{equation}
It is interesting to note that while previous studies (see \cite{MOR_H}-\cite{MOR_M}) of the finite population model lead to the scaling $\rho-1=O(N^{-1/2})$, the spectrum involves the transition scale $\rho-1=O(N^{-1/3})$. 

Our basic approach is to use singular perturbation methods to analyze the system of ODEs when $N$ becomes large. The problem can then be reduced to solving simpler, single differential equations whose solutions are known, such as Hermite equations. Our analysis does make some assumptions about the forms of various asymptotic series, and about the asymptotic matching of expansions on different scales. We also comment that we assume that the eigenvalue index $j$ is $O(1)$; thus we are not computing the large eigenvalues here. 

This paper is organized as follows. In section 2 we state the mathematical problem and obtain the basic equations. In section 3 we summarize our final asymptotic results for the eigenvalues and the (unnormalized) eigenvectors. The derivations are relegated to section 4. Some numerical studies appear in section 5; these assess the accuracy of the asymptotics.

\section{Statement of the problem}
Throughout the paper we assume that time has been scaled so to make the mean service time $1/\mu=1$.

The conditional sojourn time density $p_n(t)$ satisfies the following linear system of ordinary differential equations, or, equivalently, differential-difference equation:
\begin{equation}\label{s2_fs_rec}
p_n'(t)=\rho\Big(1-\frac{n+1}{N}\Big)p_{n+1}(t)+\frac{n}{n+1}p_{n-1}(t)-\Big[\rho\Big(1-\frac{n+1}{N}\Big)+1\Big]p_n(t),\;0\le n\le N-1.
\end{equation}
The above holds also at $n=0$ if we require $p_{-1}$ to be finite, and when $n=N-1$ (\ref{s2_fs_rec}) becomes 
\begin{equation}\label{s2_fs_p-1}
p_{N-1}'(t)=\frac{N-1}{N}p_{N-2}(t)-p_{N-1}(t).
\end{equation}
The initial condition at $t=0$ is 
\begin{equation}\label{s2_fs_pn}
p_n(0)=\frac{1}{n+1},
\end{equation}
and we note that (\ref{s2_fs_pn}) follows from integrating (\ref{s2_fs_rec}) from $t=0$ to $t=\infty$.

Here we focus on obtaining the eigenvalues (and corresponding eigenvectors) for the matrix $\mathbf{A}=\mathbf{A}(N;\rho)$, which corresponds to the difference operator in the right-hand side of (\ref{s2_fs_rec}), specifically
\begin{displaymath}
\mathbf{A}=\left[\begin{array}{cccccccc}
-1-\rho\big(\frac{N-1}{N}\big) & \rho\big(\frac{N-1}{N}\big) & 0 & 0 & 0 & \cdots & 0 & 0\\
1/2 & -1-\rho\big(\frac{N-2}{N}\big) & \rho\big(\frac{N-2}{N}\big) & 0 & 0 & \cdots & 0 & 0 \\
0 & 2/3 & -1-\rho\big(\frac{N-3}{N}\big) & \rho\big(\frac{N-3}{N}\big) & 0 & \cdots & 0 & 0 \\
\vdots  & \vdots  & \vdots  & \vdots  & \vdots & \ddots & \vdots  & \vdots \\
0 & 0 & 0 &0 & 0 & \cdots & \frac{N-1}{N}  & -1
\end{array}\right]
\end{displaymath}
which has size $N\times N$. It follows that the solution of (\ref{s2_fs_rec})-(\ref{s2_fs_pn}) has the spectral representation
\begin{equation}\label{s2_fs_sum}
p_n(t)=\sum_{j=0}^{N-1} e^{-\nu_j(N,\rho) t}c_j\phi_j(n;N,\rho).
\end{equation}
Here $-\nu_j$ are the eigenvalues of $\mathbf{A}$, indexed by $0\le j\le N-1$ and ordered as $0<\nu_0<\nu_1<\nu_2<\cdots<\nu_{N-1}$, $\phi_j(n)$ is the eigenvector corresponding to eigenvalue $\nu_j$, and the spectral coefficients $c_j$ in (\ref{s2_fs_sum}) can be calculated from (\ref{s2_fs_pn}), hence
\begin{equation*}\label{s2_fs_phi}
\frac{1}{n+1}=\sum_{j=0}^{N-1}c_j\phi_j(n),\; 0\le n\le N-1.
\end{equation*}
From (\ref{s2_fs_rec}) we can easily establish orthogonality relations between the $\phi_j$, and these lead to an explicit expression for the $c_j$:
$$c_j=\frac{\sum_{n=0}^{N-1}\Big(\frac{\rho}{N}\Big)^n\frac{N!}{(N-n-1)!}\,\phi_j(n)}{\sum_{n=0}^{N-1}\Big(\frac{\rho}{N}\Big)^n(n+1)\frac{N!}{(N-n-1)!}\,\phi_j^2(n)}.$$

We note that the tail of the sojourn time, in view of (\ref{s2_fs_sum}), is given by 
\begin{equation}\label{s2_fs_sim}
p_n(t)\sim c_0\phi_0(n)e^{-\nu_0t},\; t\to\infty.
\end{equation}
This asymptotic relation holds for $n=O(1)$ and large times. Our analysis will assume that $N\gg 1$ but for sufficiently large $t$ (\ref{s2_fs_sim}) must still hold. We shall study the behavior of $\nu_j=\nu_j(N,\rho)$ for $N\to\infty$ and for various ranges of the parameter $\rho$.

We shall show that the behavior of the eigenvalues is very different for the cases $\rho<1$, $\rho>1$, and $\rho\sim 1$ (more precisely $\rho-1=O(N^{-1/3})$). Furthermore, within each range of $\rho$ the form of the expansions of the eigenfunctions $\phi_j(n)$ is different in several ranges of $n$. Our analysis will cover each of these ranges, but we restrict ourselves to the eigenvalue index $j$ being $O(1)$. Note that the matrix $\mathbf{A}$ has $N-1$ eigenvalues so that $j$ could be scaled as large as $O(N)$. We are thus calculating (asymptotically for $N\to\infty$) only the first few eigenvalues and their eigenfunctions. Obtaining, say, the large eigenvalues, would likely need a very different asymptotic analysis. When $j=O(1)$ we shall see that the eigenfunctions and their zeros are concentrated in a narrow range of $\xi=n/N$, which represents the fraction of the customer population that is using the processor. Since $\phi_j(n)$ are functions of the discrete variable $n$, by ``zeros'' we refer to sign changes of the eigenvectors. If the eigenvalue index $j$ were scaled to be also large with $N$, we would expect that these sign changes would be more frequent, and would occur throughout the entire interval $\xi\in(0,1)$. 

We comment that to understand fully the asymptotic structure of $p_n(t)=p_n(t;N,\rho)$ for $N\to\infty$ requires a much more complete analysis than simply knowing the eigenvalues/eigenfunctions, as the spectral expansion may not be useful in certain space/time ranges. However, (\ref{s2_fs_sim}) will always apply for sufficiently large times, no matter how large $N$ is, as long as $N$ is finite.

We also comment that in (\ref{s1_pt}), by results in \cite{SEV}, we have 
\begin{equation}\label{s2_fs_probN}
\Pr\big[\mathcal{N}(0^-)=n\big]=\frac{(\frac{\rho}{N})^n(N-1)!}{(N-1-n)!}\bigg/\sum_{l=0}^{N-1}\frac{(\frac{\rho}{N})^l(N-1)!}{(N-1-l)!},
\end{equation}
which says that the distribution of $\mathcal{N}(0^-)$ coincides with the steady state distribution of $\mathcal{N}(t)$ in a finite population queue with population $N-1$. The tail of the unconditional sojourn time is then 
\begin{equation}\label{s2_fs_ptsim}
p(t)\sim c_0'\bigg[\sum_{n=0}^{N-1}\phi_0(n)\frac{(\frac{\rho}{N})^n(N-1)!}{(N-1-n)!}\bigg]e^{-\nu_0t},\;t\to\infty,
\end{equation}
where $c_0'=c_0\big[\sum_{l=0}^{N-1}(\rho/N)^l(N-1)!/(N-1-l)!\big]^{-1}$.

\section{Summary of results}

We give results for the three cases: $\rho<1$, $\rho>1$ and $\rho-1=O(N^{-1/3})$. Within each case of $\rho$ the eigenvectors $\phi_j(n)$ have different behaviors in different ranges of $n$.

\subsection{The case $\rho<1$}
For $\rho<1$ the eigenvalues are given by 
\begin{equation}\label{s31_nu}
\nu_j=(1-\sqrt{\rho})^2+\frac{c_1}{\sqrt{N}}+\frac{c_2(j)}{N^{3/4}}+\frac{c_4(j)}{N}+O(N^{-5/4}),\; j\ge 0
\end{equation}
where 
\begin{equation}\label{s31_c1c2}
c_1=2\sqrt{\rho}\sqrt{1-\sqrt{\rho}},\quad c_2(j)=(2j+1)\sqrt{\rho}(1-\sqrt{\rho})^{3/4}
\end{equation}
and
\begin{equation}\label{s31_c4}
c_4(j)=\frac{\sqrt{\rho}(22\sqrt{\rho}-3\rho-15)}{16(1-\sqrt{\rho})}-\frac{3}{8}\sqrt{\rho}(1-\sqrt{\rho})j(j+1).
\end{equation}
We observe that the leading term in (\ref{s31_nu}) ($=(1-\sqrt{\rho})^2$) is independent of the eigenvalue index $j$ and corresponds to the relaxation rate in the standard $M/M/1$ queue. For the standard $M/M/1$-PS model (with an \underline{infinite} customer population), it is well known (see \cite{POL}, \cite{COH}) that the tail of the sojourn time density is 
\begin{equation}\label{s31_ptsim}
p(t)\sim ke^{-(1-\sqrt{\rho})^2t}\exp\Big[-3\Big(\frac{\pi}{2}\Big)^{2/3}\rho^{1/6}t^{1/3}\Big]t^{-5/6},\; t\to\infty
\end{equation}
where $k$ is a constant. This problem corresponds to solving an infinite system of ODEs, which may be obtained by letting $N\to\infty$ in the matrix $\mathbf{A}$ (with a fixed $\rho<1$). The spectrum of the resulting infinite matrix is purely continuous, which leads to the sub-exponential and algebraic factors in (\ref{s31_ptsim}). 

The result in (\ref{s31_nu}) shows that the (necessarily discrete) spectrum of $\mathbf{A}=\mathbf{A}(N,\rho)$ has, for $\rho<1$ and $N\to\infty$, all of the eigenvalues approaching $(1-\sqrt{\rho})^2$, with the deviation from this limit appearing only in the third ($O(N^{-3/4})$) and fourth ($O(N^{-1})$) terms in the expansion in (\ref{s31_nu}). Note that $c_1$ in (\ref{s31_nu}) is independent of $j$. Comparing (\ref{s31_ptsim}) to (\ref{s2_fs_ptsim}) with (\ref{s31_nu}) and $j=0$, we see that the factors $\exp\big[-3(\pi/2)^{2/3}\rho^{1/6}t^{1/3}\Big]t^{-5/6}$ in (\ref{s31_ptsim}) are replaced by $\exp\big[\big(c_1N^{-1/2}+c_2(0)N^{-3/4}+c_4(0)N^{-1}\big)t\big]$. Note that (\ref{s31_ptsim}) corresponds to letting $N\to\infty$ and then $t\to\infty$ in the finite population model, while (\ref{s2_fs_ptsim}) has $t\to\infty$ with a finite large $N$. The expansion in (\ref{s31_nu}) breaks down when $j$ becomes very large, and we note that when $j=O(N^{1/4})$, the three terms $c_1N^{-1/2}$, $c_2N^{-3/4}$ and $c_4N^{-1}$ become comparable in magnitude. The expansion in (\ref{s31_nu}) suggests that $j$ can be allowed to be slightly large with $N$, but certainly not as large as $O(N)$, which would be needed to calculate all of the eigenvalues of $\mathbf{A}$. As we stated before, here we do not attempt to get the eigenvalues of large index $j$. 

Now consider the eigenvectors $\phi_j(n)=\phi_j(n;N,\rho)$, for $\rho<1$ and $N\to\infty$. These have the expansion
\begin{equation}\label{s31_phi}
\phi_j(n)=k_0\rho^{-n/2}\Big[\Phi_j^{(0)}(y)+N^{-1/8}\Phi_j^{(1)}(y)+O(N^{-1/4})\Big]
\end{equation}
where $n$ and $y$ are related by 
\begin{equation}\label{s31_ny}
n=\frac{\sqrt{N}}{\sqrt{1-\sqrt{\rho}}}+N^{3/8}y
\end{equation}
and
\begin{equation}\label{s31_Phi0}
\Phi_j^{(0)}(y)=e^{-z^2/4}\mathrm{He}_j(z),\; z=\sqrt{2}(1-\sqrt{\rho})^{3/8}y,
\end{equation}
\begin{equation}\label{s31_Phi1}
\Phi_j^{(1)}(y)=\Big[\frac{2}{3}\beta z^2+\frac{2}{3}\beta (8j+4)\Big]\frac{d}{dz}\Big[e^{-z^2/4}\mathrm{He}_j(z)\Big]+\Big(\frac{\alpha}{2}-\frac{2}{3}\beta\Big)ze^{-z^2/4}\mathrm{He}_j(z),
\end{equation}
with
\begin{equation}\label{s31_ab}
\alpha=\sqrt{\frac{\rho}{2}}(1-\sqrt{\rho})^{-7/8},\quad \beta=-\frac{(1-\sqrt{\rho})^{1/8}}{4\sqrt{2}}.
\end{equation}
Here $\mathrm{He}_j(\cdot)$ is the $j^{th}$ Hermite polynomial, so that $\mathrm{He}_0(z)=1$ and $\mathrm{He}_1(z)=z$. The constant $k_0=k_0(j;N,\rho)$ is a normalization constant which depends upon $j$, $\rho$ and $N$, but not $n$. Apart from the factor $\rho^{-n/2}$ in (\ref{s31_phi}) we see that the eigenvectors are concentrated in the range $z,\,y=O(1)$ and this corresponds to $n=\sqrt{N}/\sqrt{1-\sqrt{\rho}}+O(N^{3/8})$. The zeros of the Hermite polynomials correspond to sign changes (with $n$) of the eigenvectors, and these are thus spaced $O(N^{3/8})$ apart. 

We next give expansions of $\phi_j(n)$ on other spatial scales, such as $n=O(1)$, $n=O(\sqrt{N})$ and $n=O(N)$, where (\ref{s31_phi}) ceases to be valid. These results will involve normalization constants that we denote by $k_l$, but each of these will be related to $k_0$, so that our results are unique up to a \underline{single} multiplicative constant. Note that the eigenvalues of $\mathbf{A}$ are all simple, which can be shown by a standard Sturm-Liouville type argument. 

For $n=\sqrt{N}x=O(\sqrt{N})$ (hence $0<x<\infty$) we find that 
\begin{equation}\label{s31_phisim}
\phi_j (n)\sim k_1\rho^{-n/2}g_j(x)e^{N^{1/4}f(x)},\; N\to\infty,\; 0<x=\frac{n}{\sqrt{N}}<\infty
\end{equation} 
where 
\begin{equation*}
f(x)=2\sqrt{x}-\frac{2}{3}\sqrt{1-\sqrt{\rho}}\,x^{3/2},
\end{equation*}
\begin{equation}\label{s31_fg}
g_j(x)=\frac{(1-\sqrt{1-\sqrt{\rho}}\,x)^je^{x^2/4}}{x^{1/4}[1+(1-\sqrt{\rho})^{1/4}\sqrt{x}]^{2j+1}}\exp\Big[(2j+1)(1-\sqrt{\rho})^{1/4}\sqrt{x}\Big].
\end{equation}
The normalization constants $k_0$ and $k_1$ are related by 
\begin{equation}\label{s31_k0}
k_0= k_1N^{-j/8}(-1)^j2^{-5j/2-1}(1-\sqrt{\rho})^{{(j+1)}/{8}}\exp\Big[\frac{1}{4(1-\sqrt{\rho})}+2j+1+\frac{4}{3}(1-\sqrt{\rho})^{-1/4}N^{1/4}\Big],
\end{equation}
by asymptotic matching between (\ref{s31_phi}) and (\ref{s31_phisim}).

Next we consider $n=O(N)$ and scale $n=N\xi$. The leading term becomes
\begin{equation}\label{s31_phileading}
\phi_j(n)\sim k_2\rho^{-n/2}G(\xi,j)\exp\big[NF(\xi)+N^{1/2}F_1(\xi)+N^{1/4}F_2(\xi,j)\big]
\end{equation}
with
\begin{eqnarray}\label{s31_F}
F(\xi)&=&\frac{1}{2}\xi\log\rho+\xi\log 2+\frac{1}{2}\xi-\frac{1}{2\sqrt{\rho}}\sqrt{\rho\xi^2+4\xi(1-\sqrt{\rho})}\\
&&-\xi\log\Big[2\sqrt{\rho}-\rho\xi+\sqrt{\rho^2\xi^2+4\rho\xi(1-\sqrt{\rho})}\Big]\nonumber\\
&&+\frac{\rho-4\sqrt{\rho}+2}{2\rho}\log\Big[\frac{\rho\xi+2-2\sqrt{\rho}+\sqrt{\rho^2\xi^2+4\rho\xi(1-\sqrt{\rho})}}{2(1-\sqrt{\rho})}\Big]\nonumber\\
&&+\frac{1}{2}\log\bigg[\frac{(\rho-2\sqrt{\rho}+2)\xi+2(1-\sqrt{\rho})-(\sqrt{\rho}-2)\sqrt{\rho\xi^2+4\xi(1-\sqrt{\rho})}}{2(1-\sqrt{\rho})}\bigg],\nonumber
\end{eqnarray}
\begin{equation}\label{s31_F1}
F_1(\xi)=\frac{c_1}{\rho}\log\bigg[\frac{2(1-\sqrt{\rho})+\rho\xi+\sqrt{\rho^2\xi^2+4\rho\xi(1-\sqrt{\rho})}}{2(1-\sqrt{\rho})}\bigg],
\end{equation}
\begin{equation}\label{s31_F2}
F_2(\xi,j)=\frac{c_2(j)}{\rho}\log\bigg[\frac{2(1-\sqrt{\rho})+\rho\xi+\sqrt{\rho^2\xi^2+4\rho\xi(1-\sqrt{\rho})}}{2(1-\sqrt{\rho})}\bigg]
\end{equation}
and $G(\xi,j)$ can be written as the integral
\begin{equation}\label{s31_G}
G(\xi,j)=\xi^{-3/4}\exp\Big[\int_0^\xi\Big(\frac{3}{4v}-H(v,j)\Big)dv\Big]
\end{equation}
where 
\begin{equation*}\label{s31_H}
H(\xi,j)=\frac{1}{(1-\xi)e^{F'}-e^{-F'}}\bigg\{\frac{c_4(j)+\rho}{\sqrt{\rho}}-e^{F'}-\frac{1}{\xi}e^{-F'}+\Big(1-\frac{\sqrt{\rho}}{2}\xi\Big)\big[F''+(F_1')^2\big]\bigg\}
\end{equation*}
and $F'$ and $F_1'$ are, respectively, the derivatives of (\ref{s31_F}) and (\ref{s31_F1}). Note that $F$ and $F_1$ are independent of the eigenvalue index $j$, while $F_2$ and $G$ do depend on it. The constants $k_2$ and $k_1$ are related by  
%\begin{equation*}\label{s31_k2}
$k_2=(-1)^j(1-\sqrt{\rho})^{-1/4}N^{-3/8}k_1$,
%\end{equation*} 
by asymptotic matching between (\ref{s31_phisim}) and (\ref{s31_phileading}).

Finally we consider the scale $n=O(1)$. The expansion in (\ref{s31_phisim}) with (\ref{s31_fg}) develops a singularity as $x\to 0^+$ and ceases to be valid for small $x$. For $n=O(1)$ we obtain
\begin{equation}\label{s31_phisim3}
\phi_j(n)\sim k_3\rho^{-n/2}\frac{1}{2\pi i}\oint\frac{1}{z^{n+1}}\frac{1}{1-z}\exp\Big(\frac{1}{1-z}\Big)dz,
\end{equation}
where the contour integral is a small loop about $z=0$, and by asymptotic matching of (\ref{s31_phisim}) as $x\to 0$ with (\ref{s31_phisim3}) as $n\to \infty$, $k_3=k_1{2\sqrt{\pi}}N^{1/8}/{\sqrt{e}}.$

\subsection{The case $\rho>1$}
We next consider $N\to\infty$ with $\rho>1$. The eigenvalues are now small, of the order $O(N^{-1})$, with
\begin{equation}\label{s31_nuj}
\nu_j=\frac{1}{N}\Big(\frac{\rho}{\rho-1}+\rho j\Big)+o(N^{-1}).
\end{equation}
Note that now we do not see the coalescence of eigenvalues, as was the case when $\rho<1$, and the eigenvalue index $j$ appears in the leading term in (\ref{s31_nuj}). The form in (\ref{s31_nuj}) also suggests that the tail behavior in (\ref{s2_fs_sim}) is achieved when $t/N\gg 1$. Now the zeros of the eigenvectors will be concentrated in the range where $n=N(1-\rho^{-1})+O(\sqrt{N})$, and introducing the new spatial variable $X$, with
\begin{equation}\label{s31_X}
n=N\Big(1-\frac{1}{\rho}\Big)+\sqrt{N}X,
\end{equation}
we find that 
\begin{equation}\label{s31_phijX}
\phi_j(n)\sim k_0\mathrm{He}_j(\sqrt{\rho} X)
\end{equation}
where $\mathrm{He}_j(\cdot)$ is again the $j^{th}$ Hermite polynomial, and $k_0$ is again a normalization constant, possibly different from (\ref{s31_phi}). On the $\xi$-scale with $\xi=n/N$ we obtain 
\begin{equation}\label{s31_phijk1}
\phi_j(n)\sim k_1\xi^{\frac{1}{\rho-1}}\Big(\xi-1+\frac{1}{\rho}\Big)^j,\; \xi\ne 1-\frac{1}{\rho}
\end{equation}
and $k_0$ and $k_1$ are related by
\begin{equation}\label{s31_k0phi}
k_0=\rho^{-j/2}\Big(1-\frac{1}{\rho}\Big)^{\frac{1}{\rho-1}}N^{-j/2}k_1,
\end{equation}
by asymptotic matching between (\ref{s31_phijX}) and (\ref{s31_phijk1}).
Note that $\xi=1-1/\rho+X/\sqrt{N}$ and for the first two eigenvectors ($j=0,\,1$), (\ref{s31_phijX}) is a special case of (\ref{s31_phijk1}). The expression in (\ref{s31_phijk1}) holds for $0<\xi<1-1/\rho$ and for $1-1/\rho<\xi< 1$, but not for $n=O(1)\;(\xi=O(N^{-1}))$ or $\xi\sim 1-1/\rho$. For the latter we must use (\ref{s31_phijX}) and for $n=O(1)$ we shall show that 
\begin{equation}\label{s31_phiint}
\phi_j(n)\sim k_2\frac{e^{i\pi\rho/(\rho-1)}}{2\pi i}\int_\mathcal{C}(1-z)^{-\frac{\rho}{\rho-1}}\Big(z-\frac{1}{\rho}\Big)^{\frac{1}{\rho-1}}z^ndz,
\end{equation}
where $\mathcal{C}$ is a closed loop that encircles the branch cut, where $\Im(z)=0$ and $\Re(z)\in[\rho^{-1},1]$, in the $z$-plane, with the integrand being analytic exterior to this cut. By expanding (\ref{s31_phiint}) for $n\to\infty$ and matching to (\ref{s31_phijk1}) as $\xi\to 0^+$ we obtain
\begin{equation}\label{s31_k1k2}
k_2=k_1N^{\frac{1}{1-\rho}}(-1)^j\Big(1-\frac{1}{\rho}\Big)^{j-\frac{1}{\rho-1}}\Gamma\Big(\frac{\rho}{\rho-1}\Big).
\end{equation}

In contrast to when $\rho<1$, the eigenvectors now vary smoothly throughout the interval $\xi\in(0,1)$ (cf. (\ref{s31_phijk1})) but their sign changes are all concentrated where $\xi\sim 1-1/\rho$ and the spacings of these changes are of the order $O(N^{1/2})$, and approximately the same as the spacings of the zeros of the Hermite polynomials (cf. (\ref{s31_X}) and (\ref{s31_phijX})).

\subsection{The case $\rho\approx 1$}
Finally we consider the case $\rho\sim 1$ and introduce the parameter $\gamma$ by 
\begin{equation*}\label{s31_gamma}
\rho=1+\frac{\gamma}{N^{1/3}},\;-\infty<\gamma<\infty.
\end{equation*}
This case will asymptotically match, as $\gamma\to +\infty$, to the $\rho>1$ results and, as $\gamma\to -\infty$, to the $\rho<1$ results. Now a two-term asymptotic approximation to the eigenvalues is 
\begin{equation}\label{s31_nu2term}
\nu_j=f(\gamma)N^{-2/3}+g(j,\gamma)N^{-1}+o(N^{-1})
\end{equation} 
where
\begin{equation}\label{s31_gammaA}
\gamma-A+\frac{2}{A^2}=0,
\end{equation}
\begin{equation}\label{s31_fA}
f(\gamma)=\frac{1}{A}+\frac{1}{A^4},
\end{equation}
\begin{equation}\label{s31_gA}
g(j,\gamma)=\frac{1}{A^6}-\frac{3}{A^3}+\frac{1}{2}+\Big(j+\frac{1}{2}\Big)\sqrt{1+\frac{4}{A^3}}.
\end{equation}
Thus given $\gamma$ we must solve (\ref{s31_gammaA}) to get $A=A(\gamma)$ and then compute $f$ and $g$. We can explicitly invert (\ref{s31_gammaA}) to obtain
\begin{equation}\label{s31_Ab}
A(\gamma)=\frac{\gamma^2+b\gamma+b^2}{3b}\quad\textrm{with}\quad b=\big[\gamma^3+27+3\sqrt{81+6\gamma^3}\big]^{1/3}.
\end{equation}
We also note that if $\rho=1$ we have $\gamma=0$ and then $A=A(0)=2^{1/3}$, $f(0)=3\cdot 2^{-4/3}$, and $g(j,0)=-3/4+(j+1/2)\sqrt{3}$. 
The expression in (\ref{s31_nu2term}) shows that the eigenvalues are small, of the order $O(N^{-2/3})$, and to leading order coalesce at $f(\gamma)N^{-2/3}$. The second term, however, depends linearly on the eigenvalue index $j$. 

The expansions of the eigenvectors will now be different on the four scales $n=O(1)$, $n=O(N^{2/3})$, $n=N^{2/3}A(\gamma)+O(\sqrt{N})$ and $n=O(N)$. It is on the third scale that the zeros of $\phi_j(n)$ become apparent, and if we introduce $U$ by 
\begin{eqnarray}\label{s31_U}
n=N^{2/3}A(\gamma)+\sqrt{N}U,\; -\infty<U<\infty
\end{eqnarray}
we obtain the following leading order approximation to the eigenvectors 
\begin{equation}\label{s31_phieig}
\phi_j(n)\sim k_0\exp\bigg\{\frac{U}{A^2}N^{1/6}+\frac{U^2}{4}\Big(1-\sqrt{1+\frac{4}{A^3}}\Big)\bigg\}\mathrm{He}_j\Big(\Big(1+\frac{4}{A^3}\Big)^{1/4}U\Big).
\end{equation}
Thus again the Hermite polynomials arise, but now on the scale $U=O(1)$, which corresponds to $n-N^{2/3}A(\gamma)=O(\sqrt{N})$. The spacing of the zeros (or sign changes) of $\phi_j(n)$ for the case $\rho-1=O(N^{-1/3})$ is thus $O(\sqrt{N})$ which is comparable to the case $\rho>1$, and unlike the case $\rho<1$ where (cf. (\ref{s31_ny})) the spacing was $O(N^{3/8})$. 

On the spatial scale $n=N^{2/3}V=O(N^{2/3})$, (\ref{s31_phieig}) ceases to be valid and then we obtain 
\begin{equation}\label{s31_phiV}
\phi_j(n)\sim k_1\mathcal{G}(V,j)e^{N^{1/3}\mathcal{F}(V)},\; V=\frac{n}{N^{2/3}}\ne A(\gamma)
\end{equation}
where this applies for all $V>0$ except $V\sim A(\gamma)$, where (\ref{s31_phieig}) holds, and $\mathcal{F}$ is given by 
\begin{eqnarray}\label{s31_mathcalF}
\mathcal{F}(V)=\left\{ \begin{array}{ll}
\int_0^V\mathcal{F}'_+(w)dw, & 0\le V\le A(\gamma)\\
\int_0^A\mathcal{F}'_+(w)dw+\int_A^V\mathcal{F}'_-(w)dw, & V> A(\gamma),
\end{array} \right.
\end{eqnarray}
\begin{equation}\label{s31_F'pm}
\mathcal{F}'_\pm(V)=\frac{1}{2}\Big[V-\gamma\pm\sqrt{(V-\gamma)^2-4\big(f(\gamma)-\frac{1}{V}\big)}\Big],
\end{equation}
where $f(\gamma)=A^{-1}+A^{-4}$ is as in (\ref{s31_fA}). Note that $\mathcal{F}$ is independent of the eigenvalue index $j$, and if $V=A(\gamma)$ then $(A-\gamma)^2=4(f(\gamma)-A^{-1})$, which follows from (\ref{s31_gammaA}) and (\ref{s31_fA}). Thus the discriminant in (\ref{s31_F'pm}) vanishes when $V=A(\gamma)$, and in fact it has a double zero at this point. Thus the sign switch in $\mathcal{F}'$ as $V$ crosses $A(\gamma)$ is needed to smoothly continue this function from $V<A$ to $V>A$. The function $\mathcal{G}(V, j)$ is given by
\begin{equation}\label{s31_Gv<a}
\mathcal{G}(V,j)=V^{-1/4}\big[V-A(\gamma)\big]^{j}\exp\bigg\{\int_0^V\Big[\frac{1}{4w}-\frac{j}{w-A}-H^{(1)}_+(w,j)\Big]dw\bigg\},\; 0\le V< A
\end{equation}
and
\begin{eqnarray}\label{s31_Gv>a}
\mathcal{G}(V,j)&=&A(\gamma)^{-1/4}\big[V-A(\gamma)\big]^{j}\exp\bigg\{\int_0^A\Big[\frac{1}{4w}-\frac{j}{w-A}-H^{(1)}_+(w,j)\Big]dw\bigg\}\nonumber\\
&&\times\exp\bigg\{-\int_A^V\Big[\frac{j}{w-A}+H^{(1)}_-(w,j)\Big]dw\bigg\},\; V> A.
\end{eqnarray}
We can show that $H^{(1)}_+(V,j)\sim(4V)^{-1}$ as $V\to 0$ and $H^{(1)}_\pm(V,j)\sim -j/(V-A)$ as $V\to A$, so that the integrals in (\ref{s31_Gv<a}) and (\ref{s31_Gv>a}) are convergent at $w=0$ and $w=A$. 

By expanding (\ref{s31_phiV})  for $V\to A$ we can establish the asymptotic matching of the results on the $V$-scale (cf. (\ref{s31_phiV})) and the $U$-scale (cf. (\ref{s31_phieig})) and then relate the constants $k_0$ and $k_1$, leading to 
\begin{equation*}\label{s31_k0andk1}
k_0=k_1A^{-1/4}(1+4A^{-3})^{-j/4}N^{-j/6}\exp\bigg\{N^{1/3}\int_0^A\mathcal{F}_+'(w)dw+\int_0^A\Big[\frac{1}{4w}-\frac{j}{w-A}-H^{(1)}_+(w,j)\Big]dw\bigg\}.
\end{equation*}

Next we consider the scale $\xi=n/N=O(1)$ with $0<\xi\le 1$, and now the eigenvectors have the expansion
\begin{equation}\label{s31_phixi}
\phi_j(n)\sim k_2G(\xi,j)e^{N^{1/3}F(\xi)}
\end{equation}
where
\begin{equation}\label{s31_Fxi}
F(\xi)=f(\gamma)\log \xi=\Big(\frac{1}{A}+\frac{1}{A^4}\Big)\log\xi
\end{equation}
and 
\begin{equation}\label{s31_Gxi}
G(\xi,j)=\xi^{g(j,\gamma)-\gamma f(\gamma)}\exp\Big[\frac{1-\gamma f(\gamma)}{\xi}\Big]
\end{equation}
where $g$ and $f$ are as in (\ref{s31_fA}) and (\ref{s31_gA}). The constants $k_2$ and $k_1$ are related by
\begin{eqnarray*}\label{s31_k2andk1}
k_2&=&k_1N^{[g(j,\gamma)-\gamma f(\gamma)]/3}A^{-g(j,\gamma)+\gamma f(\gamma)+j-1/4}\exp\bigg\{\int_0^A\Big[\frac{1}{4w}-\frac{j}{w-A}-H^{(1)}_+(w,j)\Big]dw\bigg\}\nonumber\\
&&\times\exp\bigg\{N^{1/3}\Big[\frac{f(\gamma)}{3}\log N-f(\gamma)\log A-\frac{1-\gamma f(\gamma)}{A}+\int_0^A\mathcal{F}_+'(w)dw\Big]\bigg\}.
\end{eqnarray*}

For $n=O(1)$, (\ref{s31_phisim3}) holds with $\rho^{-n/2}\sim 1$, but now
$k_3=2\sqrt{\pi/e}\,N^{1/6}(-1)^j\big[A(\gamma)\big]^jk_1$.

To summarize, we have shown that the eigenvalues have very different behaviors for $\rho<1$ (cf. (\ref{s31_nu})), $\rho>1$ (cf. (\ref{s31_nuj})) and $\rho-1=O(N^{-1/3})$ (cf. (\ref{s31_nu2term})). In the first case, as $N\to\infty$, the eigenvalues all coalesce about $(1-\sqrt{\rho})^2$ which is the relaxation rate of the standard $M/M/1$ model. However, higher order terms in the expansion show the splitting of the eigenvalues (cf. $c_2(j)$ and $c_4(j)$ in (\ref{s31_nu})), which occurs at the $O(N^{-3/4})$ term in the expansion of the $\nu_j$. For $\rho>1$ the eigenvalues are small, of order $O(N^{-1})$, but even the leading term depends upon $j$ (cf. (\ref{s31_nuj})). When $\rho-1=\gamma N^{-1/3}$ there is again a coalescence of the eigenvalues, now about $N^{-2/3}(A^{-1}+A^{-4})$ where $A=A(\gamma)$ is given by (\ref{s31_Ab}). The splitting now occurs at the first correction term, which is of order $O(N^{-1})$. We also note that the leading order dependence of the $\nu_j=\nu_j(N,\rho)$ on $j$ occurs always in a simple linear fashion. Our analysis will also indicate how to compute higher order terms in the expansions of the $\nu_j$ and $\phi_j(n)$, for all three cases of $\rho$ and all ranges of $n$. Ultimately, obtaining the leading terms for the $\phi_j(n)$ reduces in all 3 cases of $\rho$ to the classic eigenvalue problem for the quantum harmonic oscillator, which can be solved in terms of Hermite polynomials.

\section{Brief derivations	}
We proceed to compute the eigenvalues and eigenvectors of the matrix $\mathbf{A}$ above (\ref{s2_fs_sum}), treating respectively the cases $\rho<1$, $\rho>1$ and $\rho-1=O(N^{-1/3})$, in subsections 4.1-4.3. We always begin by considering the scaling of $n$, with $N$, where the oscillations or sign changes of the eigenvectors occur, and this scale also determines the asymptotic eigenvalues. Then, other spatial ranges of $n$ will be treated, which correspond to the ``tails'' of the eigenvectors. 

\subsection{The case $\rho<1$}
We recall that $\rho=\lambda_0 N/\mu$ so that $\rho<1$ means that the service rate $\mu$ exceeds the maximum total arrival rate $\lambda_0 N$. When $\rho<1$ and $N\to\infty$ the distribution of $\mathcal{N}(0^-)$ in (\ref{s2_fs_probN}) behaves as $\Pr\big[\mathcal{N}(0^-)=n\big]\sim (1-\rho)\rho^n$ for $n=O(1)$, which is the same as the result for the infinite population $M/M/1$-PS queue. 

We introduce $$y=\bigg(n-\frac{\sqrt{N}}{\sqrt{1-\sqrt{\rho}}}\bigg)N^{-3/8}=nN^{-3/8}-\frac{N^{1/8}}{\sqrt{1-\sqrt{\rho}}}$$ 
as in (\ref{s31_ny}), and set 
%\begin{equation*}\label{s41_phiPhiy}
$\phi_j(n;N,\rho)=\rho^{-n/2}\Phi_j(y;N,\rho)$
%\end{equation*}
and
%\begin{equation*}\label{s41_nuy}
$\nu_j(N,\rho)=(1-\sqrt{\rho})^2+{c(N,\rho)}/{\sqrt{N}}$.
%\end{equation*}
Then setting $p_n(t)=e^{-\nu_j t}\rho^{-n/2}\Phi_j(y)$ in (\ref{s2_fs_rec}) and noting that changing $n$ to $n\pm 1$ corresponds to changing $y$ to $y\pm N^{-3/8}$, we find that 
\begin{eqnarray}\label{s41_expany}
&&-c\bigg(\frac{1}{\sqrt{1-\sqrt{\rho}}}+\frac{y}{N^{1/8}}+\frac{1}{\sqrt{N}}\bigg)\Phi_j(y)\\
&=&\sqrt{\rho}\bigg(\frac{\sqrt{N}}{\sqrt{1-\sqrt{\rho}}}+yN^{3/8}+1\bigg)\bigg(1-\frac{1}{\sqrt{N}\sqrt{1-\sqrt{\rho}}}-\frac{y}{N^{5/8}}-\frac{1}{N}\bigg)\Phi_j(y+N^{-3/8})\nonumber\\
&&+\;\sqrt{\rho}\bigg(\frac{\sqrt{N}}{\sqrt{1-\sqrt{\rho}}}+yN^{3/8}\bigg)\Phi_j(y-N^{-3/8})\nonumber\\
&&+\;\bigg[\rho\bigg(\frac{1}{\sqrt{N}\sqrt{1-\sqrt{\rho}}}+\frac{y}{N^{5/8}}+\frac{1}{N}\bigg)-2\sqrt{\rho}\bigg]\bigg(\frac{\sqrt{N}}{\sqrt{1-\sqrt{\rho}}}+yN^{3/8}+1\bigg)\Phi_j(y).\nonumber
\end{eqnarray}
Here we also multiplied (\ref{s2_fs_rec}) by $n+1$, which will simplify some of the expansions that follow.

Letting $N\to\infty$ in (\ref{s41_expany}) leads to $-c(1-\sqrt{\rho})^{-1/2}\Phi_j(y)=-2\sqrt{\rho}\Phi_j(y)+O(N^{-1/8})$ so we conclude that 
%\begin{equation*}\label{s41_ctoy}
$c(N,\rho)\to 2\sqrt{\rho}\sqrt{1-\sqrt{\rho}}$ as $N\to\infty$ and $\rho<1$.
%\end{equation*}
The form of (\ref{s41_expany}), which has the small parameter $N^{-1/8}$, then suggests that we expand $\Phi_j(y)$ as
\begin{equation}\label{s41_expandPhiy}
\Phi_j(y)=\Phi_j^{(0)}(y)+N^{-1/8}\Phi_j^{(1)}(y)+N^{-1/4}\Phi_j^{(2)}(y)+O(N^{-3/8}),
\end{equation} 
and we also expand $c=c(N,\rho)$ as 
\begin{equation}\label{s41_expandcy}
c=2\sqrt{\rho}\sqrt{1-\sqrt{\rho}}+\frac{c_2}{N^{1/4}}+\frac{c_3}{N^{3/8}}+\frac{c_4}{N^{1/2}}+O(N^{-5/8}).
\end{equation}
Then we obtain from (\ref{s41_expany}) the limiting ODE
\begin{equation}\label{s41_limitode}
a\sqrt{\rho}\frac{d^2}{dy^2}\Phi_j^{(0)}(y)+\big[ac_2-\sqrt{\rho}(1-\sqrt{\rho})y^2\big]\Phi_j^{(0)}(y)=0,
\end{equation}
where we set, for convenience, $a=1/\sqrt{1-\sqrt{\rho}}$.

Letting $\mathcal{L}$ be the differential operator
$\mathcal{L}\big\{f(y)\big\}=a\sqrt{\rho}f''(y)+\big[ac_2-\sqrt{\rho}(1-\sqrt{\rho})y^2\big]f(y)$, at the next two orders ($O(N^{-3/8})$ and $O(N^{-1/2})$) we then obtain
\begin{equation}\label{s41_order1}
\mathcal{L}\Big\{\Phi_j^{(1)}\Big\}+(yc_2+ac_3)\Phi_j^{(0)}+\sqrt{\rho}(1-a^2)\frac{d}{dy}\Phi_j^{(0)}+y\sqrt{\rho}\frac{d^2}{dy^2}\Phi_j^{(0)}=0
\end{equation}
and
\begin{equation}\label{s41_order2}
\mathcal{L}\Big\{\Phi_j^{(2)}\Big\}+(yc_2+ac_3)\Phi_j^{(1)}+(yc_3+ac_4)\Phi_j^{(0)}+\sqrt{\rho}(1-a^2)\frac{d}{dy}\Phi_j^{(1)}+\;y\sqrt{\rho}\frac{d^2}{dy^2}\Phi_j^{(1)}-2\sqrt{\rho}ay\Phi_j^{(0)}=0.
\end{equation}
Note that the coefficients $c_l$ in (\ref{s41_expandcy}) will depend upon $\rho$ and also the eigenvalue index $j$. We also require the eigenfunctions $\Phi_j(y)$ to decay as $y\to\pm\infty$. 

Changing variables from $y$ to $z$ with 
$y=(1-\sqrt{\rho})^{-3/8}z/\sqrt{2}=a^{3/4}z/\sqrt{2}$
we see that 
$\mathcal{L}\big\{f(y)\big\}=2\sqrt{\rho/a}\,L\big\{f(z)\big\}$,
where
$L\big\{f(z)\big\}=f''(z)+\big[a^{3/2}\,c_2/(2\sqrt{\rho})-z^2/4\big]f(z)$.
Thus solving (\ref{s41_limitode}) corresponds to $L\big\{\Phi_j^{(0)}\big\}=0$, which is a standard eigenvalue problem. The only acceptable solutions (which decay as $z\to\pm\infty$) correspond to
\begin{equation}\label{s41_c2}
\frac{a^{3/2}}{2\sqrt{\rho}}c_2=j+\frac{1}{2};\quad j=0,1,2,\cdots
\end{equation}
and the (unnormalized) eigenfunctions are 
%\begin{equation}\label{s41_eigenfunction}
$\Phi_j^{(0)}=e^{-z^2/4}\mathrm{He}_j(z),\; j=0,1,2,\cdots.$
%\end{equation}
Here $\mathrm{He}_j(\cdot)$ is the Hermite polynomial, which satisfies $\mathrm{He}_j(z)\sim z^j$ for $z\to\pm\infty$.

We proceed to compute the correction term $\Phi_j^{(1)}$ in (\ref{s41_expandPhiy}), and also $c_3(j)$ and $c_4(j)$ in (\ref{s41_expandcy}). In terms of $z$, (\ref{s41_order1}) becomes
\begin{equation}\label{s41_LPhi1}
L\big[\Phi_j^{(1)}\big]=-\frac{c_3}{2\sqrt{\rho}(1-\sqrt{\rho})^{3/4}}\Phi_j^{(0)}-\frac{(1-\sqrt{\rho})^{1/8}}{4\sqrt{2}}z^3\Phi_j^{(0)}+\frac{\sqrt{\rho}}{\sqrt{2}(1-\sqrt{\rho})^{7/8}}\frac{d}{dz}\Phi_j^{(0)}
\end{equation}
where $L\big[f(z)\big]=f''(z)+\big(j+1/2-{z^2}/{4}\big)f(z)$, in view of (\ref{s41_c2}). We determine $c_3$ by a solvability condition for (\ref{s41_LPhi1}). We multiply (\ref{s41_LPhi1}) by $\Phi_j^{(0)}$ and integrate from $z=-\infty$ to $z=\infty$, and use the properties of Hermite polynomials (see \cite{mag}). Then we conclude that $c_3=0$ for all $j$. Thus there is no $O(N^{-7/8})$ term in the expansion in (\ref{s31_nu}). 

To solve for $\Phi_j^{(1)}$ we write the right-hand side of (\ref{s41_LPhi1}) as $\alpha\frac{d}{dz}\Phi_j^{(0)}+\beta z^3\Phi_j^{(0)}$ where $\alpha$ and $\beta$ are as in (\ref{s31_ab}). Then we can construct a particular solution to (\ref{s41_LPhi1}) (with $c_3=0$) in the form $\Phi_j^{(1)}=Az^2\frac{d}{dz}\Phi_j^{(0)}+Bz\Phi_j^{(0)}+C\frac{d}{dz}\Phi_j^{(0)}$ where $A$, $B$, and $C$ are determined from 
\begin{equation}\label{s41_ABC}
\frac{3}{2}A=\beta,\quad 2A+2B=\alpha,\quad -2(2j+1)A+\frac{1}{2}C=0.
\end{equation}
Solving (\ref{s41_ABC}) leads to $\Phi_j^{(1)}$ as in (\ref{s31_Phi1}). 

To compute $c_4$, the $O(N^{-1})$ term in (\ref{s31_nu}), we use the solvability condition for the equation (\ref{s41_order2}) for $\Phi_j^{(2)}$. We omit the detailed derivation. We note that (\ref{s31_c4}) is singular in the limit $\rho\uparrow 1$, while $c_2(j)$ in (\ref{s31_c1c2}) vanishes in this limit. Also, $c_4$ is quadratic in $j$ while $c_2$ is linear in $j$, so that (\ref{s31_nu}) becomes invalid both as $\rho\uparrow 1$ and as the eigenvalue index $j$ becomes large. 

We next consider the $\phi_j(n)$ on the spatial scales $n=O(\sqrt{N})$, $n=O(N)$ and $n=O(1)$. 

For $n=\sqrt{N}x=O(\sqrt{N})$ ($0<x<\infty$), the expansion in (\ref{s31_phi}) ceases to be valid. We expand the leading term $\Phi_j^{(0)}(y)$ in (\ref{s31_phi}) as $y\to\infty$ to obtain
\begin{equation}\label{s41_phix}
\phi_j(n)\sim k_02^{j/2}(1-\sqrt{\rho})^{3j/8}y^j\rho^{-n/2}\exp\Big[-\frac{1}{2}(1-\sqrt{\rho})^{3/4}y^2\Big],
\end{equation}
which suggests that we expand the eigenvector $\phi_j(n)$ in the form (\ref{s31_phisim}) on the $x$-scale, noticing also that $y=(x-1/\sqrt{1-\sqrt{\rho}})N^{1/8}$. We set $p_n(t)=e^{-\nu_jt}\phi_j(n)$ in (\ref{s2_fs_rec}) with $\nu_j$ given by (\ref{s31_nu}) and $\phi_j(n)$ having the form in (\ref{s31_phisim}). For $N\to\infty$ we obtain the following ODEs for $f(x)$ and $g_j(x)$:
\begin{equation}\label{s41_xf(x)ode}
\big[f'(x)\big]^2+(\sqrt{\rho}-1)x-\frac{1}{x}+\frac{c_1}{\sqrt{\rho}}=0
\end{equation}
and
\begin{equation}\label{s41_xg(x)ode}
g_j'(x)+\bigg[-\frac{x}{2}+\frac{3}{4x}+\frac{\sqrt{\rho}-c_2(j)x^{3/2}}{2\sqrt{\rho}x(\sqrt{1-\sqrt{\rho}}x-1)}\bigg]g_j(x)=0.
\end{equation}
Using $c_1$ and $c_2(j)$ in (\ref{s31_c1c2}), (\ref{s41_xf(x)ode}) and (\ref{s41_xg(x)ode}) can be easily solved and the results are in (\ref{s31_fg}). We note that (\ref{s41_xf(x)ode}) can be rewritten as 
$\big[f'(x)\big]^2=\big(\sqrt{1-\sqrt{\rho}}\sqrt{x}-1/\sqrt{x}\big)^2.$
After taking the square root, we choose the solution with the negative sign since $f(x)$ should achieve a maximum at $x=1/\sqrt{1-\sqrt{\rho}}$.

Now consider the scale $n=\xi N=O(N)$. Letting $x\to\infty$ in the exponential terms $N^{1/4}f(x)+x^2/4+(2j+1)(1-\sqrt{\rho})^{1/4}\sqrt{x}$ in (\ref{s31_phisim}) (with (\ref{s31_fg})) and noticing that $x=\sqrt{N}\xi$, we conclude that the expansion on the $\xi$-scale should have the form in (\ref{s31_phileading}). Using $p_n(t)=e^{-\nu_jt}\phi_j(n)$ with (\ref{s31_nu}) and (\ref{s31_phileading}) in (\ref{s2_fs_rec}) and expanding for $N\to\infty$, we obtain the following ODEs for $F$, $F_1$, $F_2$ and $G$:
\begin{equation}\label{s41_ODEf}
(1-\xi)e^{F'}+e^{-F'}+\sqrt{\rho}\xi-2=0,
\end{equation}
\begin{equation}\label{s41_ODEf1}
\Big[(1-\xi)e^{F'}-e^{-F'}\Big]F_1'+\frac{c_1}{\sqrt{\rho}}=0,
\end{equation}
\begin{equation}\label{s41_ODEf2}
\Big[(1-\xi)e^{F'}-e^{-F'}\Big]F_2'+\frac{c_2(j)}{\sqrt{\rho}}=0,
\end{equation}
and 
\begin{equation}\label{s41_ODEg}
\Big[(1-\xi)e^{F'}-e^{-F'}\Big]G'+\bigg\{\frac{c_4(j)+\rho}{\sqrt{\rho}}-e^{F'}-\frac{1}{\xi}e^{-F'}+\Big(1-\frac{\sqrt{\rho}}{2}\xi\Big)\big[F''+(F_1')^2\big]\bigg\}G=0.
\end{equation}
Solving (\ref{s41_ODEf}) for $e^{F'}$ leads to 
\begin{equation}\label{s41_ef'}
e^{F'}=\frac{2-\sqrt{\rho}\xi-\sqrt{\rho\xi^2+4\xi(1-\sqrt{\rho})}}{2(1-\xi)}.
\end{equation}
Integrating the logarithm of (\ref{s41_ef'}) leads to (\ref{s31_F}). From (\ref{s41_ef'}) we also have 
$(1-\xi)e^{F'}-e^{-F'}=-\sqrt{\rho\xi^2+4\xi(1-\sqrt{\rho})}$, 
and then solving (\ref{s41_ODEf1})-(\ref{s41_ODEg}) leads to (\ref{s31_F1})-(\ref{s31_G}).

Finally we consider the scale $n=O(1)$. We assume the leading order approximation of $\phi_j(n)$ is 
\begin{equation}\label{s41_phin}
\phi_j(n)=k_3\rho^{-n/2}\big[\psi_j(n)+o(1)\big].
\end{equation}
Then using $p_n(t)=e^{-\nu_jt}\phi_j(n)$ in (\ref{s2_fs_rec}) with $\phi_j(n)$ in (\ref{s41_phin}) and letting $N\to\infty$, we obtain the following limiting difference equation for $\psi_j(n)$:
\begin{equation}\label{s41_depsi2}
\psi_j(n+1)+\frac{n}{n+1}\psi_j(n-1)-2\psi_j(n)=0
\end{equation}
with $\psi_j(-1)$ finite (thus (\ref{s41_depsi2}) holds for all $n\ge 0$). From (\ref{s41_depsi2}) we conclude that $\psi_j(n)$ is independent of the eigenvalue index $j$, except via a multiplicative constant. Solving (\ref{s41_depsi2}) with the help of generation functions leads to (\ref{s31_phisim3}).

\subsection{The case $\rho>1$}
When $\rho>1$ the distribution of $\mathcal{N}(0^-)$ in (\ref{s2_fs_probN}) is approximately, for $N\to\infty$, a Gaussian which is centered about $n=N(1-\rho^{-1})$ and this corresponds to the fraction of the population that is typically served by the processor. We thus begin by considering the scale $n=N(1-\rho^{-1})+O(\sqrt{N})$, setting 
\begin{equation*}
\phi_j(n)=\Psi_j(X)=\Psi_j\Big(\frac{n-N(1-\rho^{-1})}{\sqrt{N}}\Big),\quad
\nu=\frac{\nu_*}{N},\quad n=N\Big(1-\frac{1}{\rho}\Big)+\sqrt{N}X.
\end{equation*}
We thus scale the eigenvalue parameter $\nu$ to be small, of order $O(N^{-1})$, which is necessary to obtain a limiting differential equation.

Setting $p_n(t)=e^{-\nu t}\Psi(X)$ and omitting for now the dependence of $\nu$ and $\Psi$ on the index $j$, (\ref{s2_fs_rec}) becomes
\begin{eqnarray}\label{s42_PsiXexpan}
-\frac{\nu_*}{N}\Psi(X)&=&\Big(1-\rho\frac{X}{\sqrt{N}}-\frac{\rho}{N}\Big)\Big[\Psi\Big(X+\frac{1}{\sqrt{N}}\Big)-\Psi(X)\Big]+\Psi\Big(X-\frac{1}{\sqrt{N}}\Big)\nonumber\\
&&-\Psi(X)-\frac{1}{N(1-\rho^{-1})+\sqrt{N}X}\Psi\Big(X-\frac{1}{\sqrt{N}}\Big).
\end{eqnarray}
Then expanding $\Psi(X)$ and $\nu_*$ as
\begin{equation}\label{s42_nuXexpan}
\Psi(X)=\Psi^{(0)}(X)+\frac{1}{\sqrt{N}}\Psi^{(1)}(X)+O(N^{-1}),\quad
\nu_*=\nu_*^{(0)}+\frac{1}{\sqrt{N}}\nu_*^{(1)}+O(N^{-1})
\end{equation}
and multiplying (\ref{s42_PsiXexpan}) by $N$ we obtain the limiting ODE
\begin{equation}\label{s42_PsiODE}
\frac{d^2}{dX^2}\Psi^{(0)}-\rho X\frac{d}{dX}\Psi^{(0)}+\Big(\nu_*^{(0)}-\frac{\rho}{\rho-1}\Big)\Psi^{(0)}=0.
\end{equation}
This can be easily transformed to the Hermite equation by setting $X=\rho^{-1/2}\widetilde{X}$. In terms of $\widetilde{X}$, (\ref{s42_PsiODE}) becomes 
\begin{equation}\label{s42_PsiXtildODE}
\frac{d^2}{d\widetilde{X}^2}\Psi^{(0)}-\widetilde{X}\frac{d}{d\widetilde{X}}\Psi^{(0)}+\frac{1}{\rho}\Big(\nu_*^{(0)}-\frac{\rho}{\rho-1}\Big)\Psi^{(0)}=0,
\end{equation}
which is the standard Hermite equation (in contrast to the parabolic cylinder equation in (\ref{s41_limitode})). Equation (\ref{s42_PsiXtildODE}) admits polynomial solutions provided that 
\begin{equation}\label{s42_nu*poly}
\frac{1}{\rho}\Big(\nu_*^{(0)}-\frac{\rho}{\rho-1}\Big)=j=0,1,2,\cdots,
\end{equation} 
which yields the leading order eigenvalue condition, and then the solution to (\ref{s42_PsiXtildODE}) is 
\begin{equation}\label{s42_Psi0}
\Psi^{(0)}(X)=k_0\mathrm{He}_j(\sqrt{\rho}X)=k_0\mathrm{He}_j(\widetilde{X}).
\end{equation}
From (\ref{s42_nu*poly}) and (\ref{s42_nuXexpan}) we have thus derived (\ref{s31_nuj}), while (\ref{s31_phijX}) follows from (\ref{s42_Psi0}). Higher order terms can be obtained by refining the expansion of (\ref{s42_PsiXexpan}) using (\ref{s42_nuXexpan}), which will lead to inhomogeneous forms of the Hermite equation; the correction terms $\nu_*^{(l)}$ for $l\ge 1$ will follow from appropriate solvability conditions.

Next we consider (\ref{s2_fs_rec}) on a broader spatial scale, introducing $\xi=n/N$, which is essentially the fraction of the population in the system (not counting the tagged customer). Letting $p_n(t)=e^{-\nu t}\varphi(\xi)$ in (\ref{s2_fs_rec}) and scaling again $\nu=\nu_*/N$ leads to 
\begin{equation}\label{s42_phixiexpan}
-\frac{\nu_*}{N}\varphi(\xi)=\Big[\rho(1-\xi)-\frac{1}{N}\Big]\Big[\varphi\Big(\xi+\frac{1}{N}\Big)-\varphi(\xi)\Big]+\Big(1-\frac{1}{N\xi+1}\Big)\varphi\Big(\xi-\frac{1}{N}\Big)-\varphi(\xi).
\end{equation}
For $N\to\infty$ (\ref{s42_phixiexpan}) leads to the limiting differential equation
\begin{equation*}\label{s42_phixiODE}
\big[\rho(1-\xi)-1\big]\varphi'(\xi)+\Big(\nu_*-\frac{1}{\xi}\Big)\varphi(\xi)=0
\end{equation*}
with solution
\begin{equation}\label{s42_phixisolu}
\varphi(\xi)=k_1\xi^{{1}/{(\rho-1)}}\Big(\xi-1+\frac{1}{\rho}\Big)^{{\nu_*}/{\rho}-{1}/{(\rho-1)}}.
\end{equation}
We argue that if the right-hand side of (\ref{s42_phixisolu}) is to be real and finite for all $\xi\in(0,1)$ we must have $\nu_*\rho^{-1}-(\rho-1)^{-1}$ a non-negative integer, and this regains the eigenvalue condition in (\ref{s31_nuj}). Alternately, since we have already fixed $\nu_*\sim \nu_*^{(0)}$ by considering the scale $X=O(1)$, we can view (\ref{s42_phixisolu}) as simply giving the approximation to the eigenfunctions on the $\xi$-scale, as given by (\ref{s31_phijk1}). Note that for $j=0$ and $j=1$ the zeros of (\ref{s42_phixisolu}) and (\ref{s42_Psi0}) coincide, but for $j\ge 2$ the expression in (\ref{s42_phixisolu}) has a zero of order $j$ at $\xi=1-\rho^{-1}$ while (\ref{s42_Psi0}) has $j$ simple zeros at points where $\sqrt{\rho}X$ is at a zero of $\mathrm{He}_j$. 

The expression in (\ref{s42_phixisolu}) vanishes as $\xi\to 0^+$ for $\rho>1$, and we thus need another expansion for the $\phi_j(n)$ for small values of $\xi$. We re-examine (\ref{s2_fs_rec}) on the scale $n=O(1)$. Setting $\nu\sim \nu_*^{(0)}/N$ and $p_n(t)=e^{-\nu t}q_n$, (\ref{s2_fs_rec}) becomes
\begin{equation*}\label{s42_qexpan}
-\nu q_n=\rho\Big(1-\frac{n+1}{N}\Big)(q_{n+1}-q_n)+\frac{n}{n+1}q_{n-1}-q_n=O\Big(\frac{1}{N}\Big)
\end{equation*}
and if $q_n\sim q_n^{(0)}$ for $N\to\infty$ with $n=O(1)$, the leading term must satisfy
\begin{equation}\label{s42_qleading}
0=\rho(q_{n+1}^{(0)}-q_{n}^{(0)})+\frac{n}{n+1}q_{n-1}^{(0)}-q_{n}^{(0)},\quad q_{-1}^{(0)}\;\textrm{finite}.
\end{equation}
Solving (\ref{s42_qleading}) using generating functions or contour integrals leads to the formula in (\ref{s31_phiint}).

\subsection{The case $\rho\approx 1$}
To analyze the cases $\rho=1$ and $\rho\approx 1$ we first note that when $\rho\uparrow 1$ the eigenvalues for $\rho<1$, which concentrate about $(1-\sqrt{\rho})^2$, behave as $\nu_j\sim(1-\rho)^2/4$ (here first $N\to\infty$, then $\rho\uparrow 1$). The result for $\rho>1$ (cf. (\ref{s31_nuj})) leads to $\nu_j\sim N^{-1}(\rho-1)^{-1}$ as $\rho\downarrow 1$, which again shows that the eigenvalues begin to coalesce. Also, $(1-\rho)^2$ balances $N^{-1}(\rho-1)^{-1}$ when $\rho-1=O(N^{-1/3})$ and this suggests the appropriate scaling for the transition region where $\rho\approx 1$. We thus define $\gamma$ by
\begin{equation}\label{s43_gamma}
\gamma=(\rho-1)N^{1/3},\; -\infty<\gamma<\infty,\; \gamma=O(1).
\end{equation}
The behavior of $\nu_j$ for $\rho\lessgtr 1$ indicates that the eigenvalues on the transition scale coalesce about some value that is $O(N^{-2/3})$ as $N\to\infty$, so we set
\begin{equation*}\label{s43_fgamma}
f=f(\gamma)=\lim_{N\to\infty}\Big(\nu_jN^{2/3}\Big);\quad \gamma,\; j=O(1).
\end{equation*}
Asymptotic matching as $\gamma\to\infty$ to the case $\rho>1$ and as $\gamma\to -\infty$ to the case $\rho<1$ leads to the following behaviors of $f(\gamma)$
\begin{equation}\label{s43_fgammasim}
f(\gamma)\sim\frac{\gamma^2}{4},\;\gamma\to -\infty;\quad f(\gamma)\sim\frac{1}{\gamma},\; \gamma\to +\infty.
\end{equation}
This analysis suggests that the eigenvalues may be expanded in the form (\ref{s31_nu2term}), where the correction term $g(j,\gamma)N^{-1}$ can also be argued by matching higher order terms in the expansions of $\nu_j$ for $\rho\lessgtr 1$.

To argue what the appropriate scaling of $n$ should be when $\rho-1=O(N^{-1/3})$, we examine (\ref{s31_ny}) as $\rho\uparrow 1$, which becomes 
\begin{equation}\label{s43_ngamma}
n=\frac{\sqrt{N}}{\sqrt{1-\sqrt{\rho}}}+N^{3/8}\frac{z}{\sqrt{2}(1-\sqrt{\rho})^{3/8}}\sim N^{2/3}\sqrt{-\frac{2}{\gamma}}+\sqrt{N}\frac{z}{2^{1/8}(-\gamma)^{3/8}}.
\end{equation}
Similarly, as $\rho\downarrow 1$ the scaling in (\ref{s31_X}) becomes
\begin{equation}\label{s43_nX}
n=N\Big(1-\frac{1}{\rho}\Big)+\sqrt{N}X\sim N^{2/3}\gamma+\sqrt{N}X.
\end{equation}
In view of (\ref{s43_ngamma}) and (\ref{s43_nX}) we scale $n$ as in (\ref{s31_U}), where $A$ is to be determined. By asymptotic matching $A=A(\gamma)$ must behave as $A(\gamma)\sim\sqrt{-2/\gamma},\;\gamma\to -\infty$ and $A(\gamma)\sim\gamma,\;\gamma\to +\infty$. 

We then set
\begin{equation}\label{s43_pntU}
p_n(t)=e^{-\nu t}\exp\big(N^{1/6}\widetilde{a}U\big)\Phi(U)
\end{equation}
where $\widetilde{a}$ is another parameter, that will be determined to insure that $\Phi$ satisfies a limiting differential equation for $N\to\infty$. For now we suppress the dependence of $\nu$ and $\Phi$ on the eigenvalue index $j$, and use (\ref{s43_pntU}) in (\ref{s2_fs_rec}) to obtain
\begin{eqnarray}\label{s43_PhiUexpan}
-\nu\Phi(U)&=&\Phi\Big(U+\frac{1}{\sqrt{N}}\Big)+\Phi\Big(U-\frac{1}{\sqrt{N}}\Big)-2\Phi(U)\nonumber\\
&&-\Big(\frac{A}{N^{1/3}}+\frac{U}{\sqrt{N}}+\frac{1}{N}\Big)\Big(1+\frac{\gamma}{N^{1/3}}\Big)\Big[\Phi\Big(U+\frac{1}{\sqrt{N}}\Big)-\Phi(U)\Big]\nonumber\\
&&-\frac{1}{N^{2/3}A+\sqrt{N}U+1}\exp\Big(-\frac{\widetilde{a}}{N^{1/3}}\Big)\Phi\Big(U-\frac{1}{\sqrt{N}}\Big)\nonumber\\
&&+\frac{\gamma}{N^{1/3}}\Big[\Phi\Big(U+\frac{1}{\sqrt{N}}\Big)-\Phi(U)\Big]\nonumber\\
&&+\Big[\exp\Big(\frac{\widetilde{a}}{N^{1/3}}\Big)-1\Big]\Big[1+\frac{\gamma}{N^{1/3}}\Big]\Phi\Big(U+\frac{1}{\sqrt{N}}\Big)\nonumber\\
&&-\Big[\exp\Big(\frac{\widetilde{a}}{N^{1/3}}\Big)-1\Big]\Big(\frac{A}{N^{1/3}}+\frac{U}{\sqrt{N}}+\frac{1}{N}\Big)\Big(1+\frac{\gamma}{N^{1/3}}\Big)\Phi\Big(U+\frac{1}{\sqrt{N}}\Big)\nonumber\\
&&+\Big[\exp\Big(-\frac{\widetilde{a}}{N^{1/3}}\Big)-1\Big]\Phi\Big(U-\frac{1}{\sqrt{N}}\Big).
\end{eqnarray}
The equation in (\ref{s43_PhiUexpan}) is an exact transformation of (\ref{s2_fs_rec}) using the scaling in (\ref{s31_U}), (\ref{s43_gamma}) and (\ref{s43_pntU}), and we rearranged the terms in the right-hand side of (\ref{s43_PhiUexpan}) in such a way that they are easier to expand for $N\to\infty$.

Next we assume that $\nu_j$ can be expanded in the form in (\ref{s31_nu2term}). With the help of Taylor expansions, we see that the right side of (\ref{s43_PhiUexpan}) will have terms that are $O(N^{-2/3})$, $O(N^{-5/6})$, and $O(N^{-1})$, with the rest being $o(N^{-1})$, and this would balance the error term(s) in the eigenvalue expansion in (\ref{s31_nu2term}). Balancing the $O(N^{-2/3})$ terms, which includes $f$ in (\ref{s31_nu2term}), leads to 
\begin{equation}\label{s43_fPhiU}
-f\Phi(U)=\Big[\widetilde{a}^2+\widetilde{a}(\gamma-A)-\frac{1}{A}\Big]\Phi(U),
\end{equation}
the $O(N^{-5/6})$ terms lead to 
\begin{equation*}\label{s43_n-5/6}
0=(\gamma-A+2\widetilde{a})\Phi'(U)+\Big(\frac{U}{A^2}-U\widetilde{a}\Big)\Phi(U)
\end{equation*}
and the $O(N^{-1})$ terms give
\begin{equation}\label{s43_n-1term}
-g\Phi(U)=\Phi''(U)-U\Phi'(U)+\Big[\frac{\widetilde{a}^2}{2}(\gamma-A)+\frac{1}{A}\Big(\widetilde{a}-\frac{U^2}{A^2}\Big)-\gamma A\widetilde{a}\Big]\Phi(U).
\end{equation}

In view of (\ref{s43_fPhiU}) we must have 
\begin{equation}\label{s43_faA}
f+\widetilde{a}^2+\widetilde{a}(\gamma-A)-\frac{1}{A}=0
\end{equation}
which is one equation relating $f$, $A$ and $\widetilde{a}$ to the ``detuning'' parameter $\gamma$. In order to obtain the second order equation in (\ref{s43_n-1term}) as the leading term we must set
%\begin{equation*}\label{s43_leadingset}
$\gamma-A+2\widetilde{a}=0$ and $A^{-2}-\widetilde{a}=0$.
%\end{equation*}
This yields two additional relations between $\widetilde{a}$ and $A$, and if we eliminate $\widetilde{a}$ we obtain precisely the equation (\ref{s31_gammaA}) that relates $\gamma$ and $A$, and then setting $\widetilde{a}=A^{-2}$ in (\ref{s43_faA}) leads to $f=A^{-1}+\widetilde{a}^2=A^{-1}+A^{-4}$, which is (\ref{s31_fA}). Finally, using $\widetilde{a}=A^{-2}$ in (\ref{s43_n-1term}) and setting $\Phi(U)=e^{U^2/4}\widetilde{\Phi}(U)$ leads to 
\begin{equation}\label{s43_Phitilde}
\widetilde{\Phi}''(U)+\Big[g-\frac{1}{2}-\frac{1}{A^6}+\frac{3}{A^3}-\Big(\frac{1}{4}+\frac{1}{A^3}\Big)U^2\Big]\widetilde{\Phi}(U)=0.
\end{equation}

If we further scale $U$ as
\begin{equation*}\label{s43_UUbar}
U=\Big(1+\frac{4}{A^3}\Big)^{-1/4}\overline{U},
\end{equation*}
(\ref{s43_Phitilde}) becomes
\begin{equation}\label{s43_ODEPhibar}
\frac{d^2}{d\overline{U}^2}\widetilde{\Phi}+\Big[\frac{g-1/2-A^{-6}+3A^{-3}}{\sqrt{1+4A^{-3}}}-\frac{1}{4}\overline{U}^2\Big]\widetilde{\Phi}=0,
\end{equation}
which is the parabolic cylinder equation in standard form. Solutions that have appropriate decay as $U$ (or $\overline{U}$) $\to\pm\infty$ require that
\begin{equation*}\label{s43_require}
\frac{g-1/2-A^{-6}+3A^{-3}}{\sqrt{1+4A^{-3}}}=j+\frac{1}{2},\quad j=0,1,2,\cdots.
\end{equation*}
Then the solution to (\ref{s43_ODEPhibar}) is proportional to a Hermite polynomial, with $\widetilde{\Phi}=\mathrm{He}_j(\overline{U})e^{-\overline{U}^2/4}$,
so we have derived (\ref{s31_gA}) and (\ref{s31_phieig}). In view of the form of (\ref{s43_PhiUexpan}), we can expand the eigenfunctions in powers of $N^{-1/6}$ and calculate higher order correction terms, obtaining an expansion of the form
$\Phi_j^{(0)}(U)+N^{-1/6}\Phi_j^{(1)}(U)+N^{-1/3}\Phi_j^{(2)}(U)+O(N^{-1/2})$.
It is likely that the expansion of the eigenvalues $\nu_j$ involves only powers of $N^{-1/3}$.

From (\ref{s31_gammaA}) we can easily obtain
%\begin{equation*}\label{s43_A(gamma)}
$A(\gamma)=\gamma+2\gamma^{-2}+O(\gamma^{-5})$ as $\gamma\to\infty$, and
$A(\gamma)=\sqrt{-2/\gamma}-\gamma^{-2}+O(|\gamma|^{-7/2})$ as $\gamma\to -\infty$,
%\end{equation*}
which can be used to verify the matching conditions in (\ref{s43_fgammasim})-(\ref{s43_nX}). We have thus shown that analysis of the case $\rho\approx 1$ is quite intricate, but ultimately, with the appropriate scaling, we again reduced the problem to the standard eigenvalue problem for the Hermite or parabolic cylinder equations. 

Next we examine the eigenvectors $\phi_j(n)$ on the scales $n=O(N^{2/3})$, $n=O(N)$ and $n=O(1)$, since (\ref{s31_phieig}) no longer applies in these ranges. We first consider the scale $n=N^{2/3}V=O(N^{2/3})$. Letting $U=N^{1/6}(V-A(\gamma))\to\infty$ in (\ref{s31_phieig}) suggests that on the $V$-scale the eigenvector $\phi_j(n)$ is in the form (\ref{s31_phiV}). Using (\ref{s31_phiV}) in (\ref{s2_fs_rec}) and noticing that changing $n$ to $n\pm 1$ corresponds to changing $V$ to $V\pm N^{-2/3}$, we obtain the following ODEs for $\mathcal{F}(V)$ and $\mathcal{G}(V,j)$:
\begin{equation}\label{s43_Fode}
\big[\mathcal{F}'(V)\big]^2+(\gamma-V)\mathcal{F}'(V)+f(\gamma)-\frac{1}{V}=0,
\end{equation}
\begin{equation}\label{s43_Gode}
\frac{d}{dV}\mathcal{G}(V,j)+\frac{g(j,\gamma)+\mathcal{F}''(V)+\big(\frac{1}{V}-\gamma V\big)\mathcal{F}'(V)+\frac{\gamma-V}{2}\big[\mathcal{F}'(V)\big]^2}{2\mathcal{F}'(V)-V+\gamma}\mathcal{G}(V,j)=0.
\end{equation}
Solving (\ref{s43_Fode}) and (\ref{s43_Gode}) leads to the results in (\ref{s31_mathcalF})-(\ref{s31_Gv>a}).

Next we consider the scale $n=\xi N=O(N)$ with $0<\xi<1$. Letting $V=N^{1/3}\xi\to \infty$ in (\ref{s31_mathcalF}) and (\ref{s31_Gv>a}), we find that the eigenvectors have an expansion in the form in (\ref{s31_phixi}). Using (\ref{s31_phixi}) in (\ref{s2_fs_rec}) yields
%\begin{equation*}
$\xi F'(\xi)-f(\gamma)=0$
%\end{equation*}
and 
\begin{equation*}
\frac{d}{d\xi}G(\xi,j)+\Big[\frac{1-\gamma f(\gamma)}{\xi^2}+\frac{\gamma f(\gamma)-g(j,\gamma)}{\xi}\Big]G(\xi,j)=0,
\end{equation*}
which can be easily solved to give (\ref{s31_Fxi}) and (\ref{s31_Gxi}). 

Finally we consider the $n=O(1)$ scale. This is necessary since $\mathcal{G}(V,j)$ in (\ref{s31_Gv<a}) is singular as $V\to 0$, with $\mathcal{G}(V,j)\sim V^{-1/4}$. We assume that the leading order approximation to $\phi_j(n)$ is $\phi_j(n)=k_3 Q_j(n)+o(1)$ and use this in (\ref{s2_fs_rec}). Similarly to the analysis of $n=O(1)$ in subsection 4.1, we find that $Q_j(n)$ satisfies $Q_j(n+1)+\frac{n}{n+1}Q_j(n-1)-2Q_j(n)=0$, which leads to the contour integral in (\ref{s31_phisim3}) with $\rho^{-n/2}\sim 1$.

\section{Numerical studies}
We assess the accuracy of our asymptotic results, and their ability to predict qualitatively and quantitatively the true eigenvalues/eigenvectors.

Recall that $\phi_j(n)$ is the $j^{\mathrm{th}}$ eigenvector with $j\ge 0$ and $0\le n\le N-1$. In Figure~\ref{f1} we plot the exact (numerical) $\phi_1(n)$ for $n\in[0,99]$, where $N=100$ and $\rho=4>1$. Our asymptotic analysis predicts that $\phi_1(n)$ will undergo a single sign change, and on the scale $X=O(1)$, (\ref{s31_phijX}) shows that $\phi_1(n)$ should be approximately linear in $X$, with a zero at $X=0$ which corresponds to $n=75$ ($=N(1-\rho^{-1})$), in view of (\ref{s31_X}). The exact eigenvector undergoes a sign change when $n$ changes from 74 to 75, in excellent agreement with the asymptotics. While the eigenvector is approximately linear near $n=75$, Figure~\ref{f1} shows that it is not globally linear, and achieves a minimum value at $n=19$. But our analysis shows that on the $\xi=n/100$ scale, we must use (\ref{s31_phijk1}), and when $\rho=4$ and $j=1$, (\ref{s31_phijk1})
%$$\xi^{\frac{1}{\rho-1}}\Big(\xi-1+\frac{1}{\rho}\Big)^j=\xi^{1/3}\Big(\xi-\frac{3}{4}\Big)$$
achieves a minimum at $\xi=3/16$, which corresponds to $n\approx 19$. This demonstrates the necessity of treating both the $X$ and $\xi$ scales. 

\begin{figure}[h]   
  \begin{minipage}[t]{0.5\linewidth}  
    \centering   
    \includegraphics[width=0.9\textwidth]{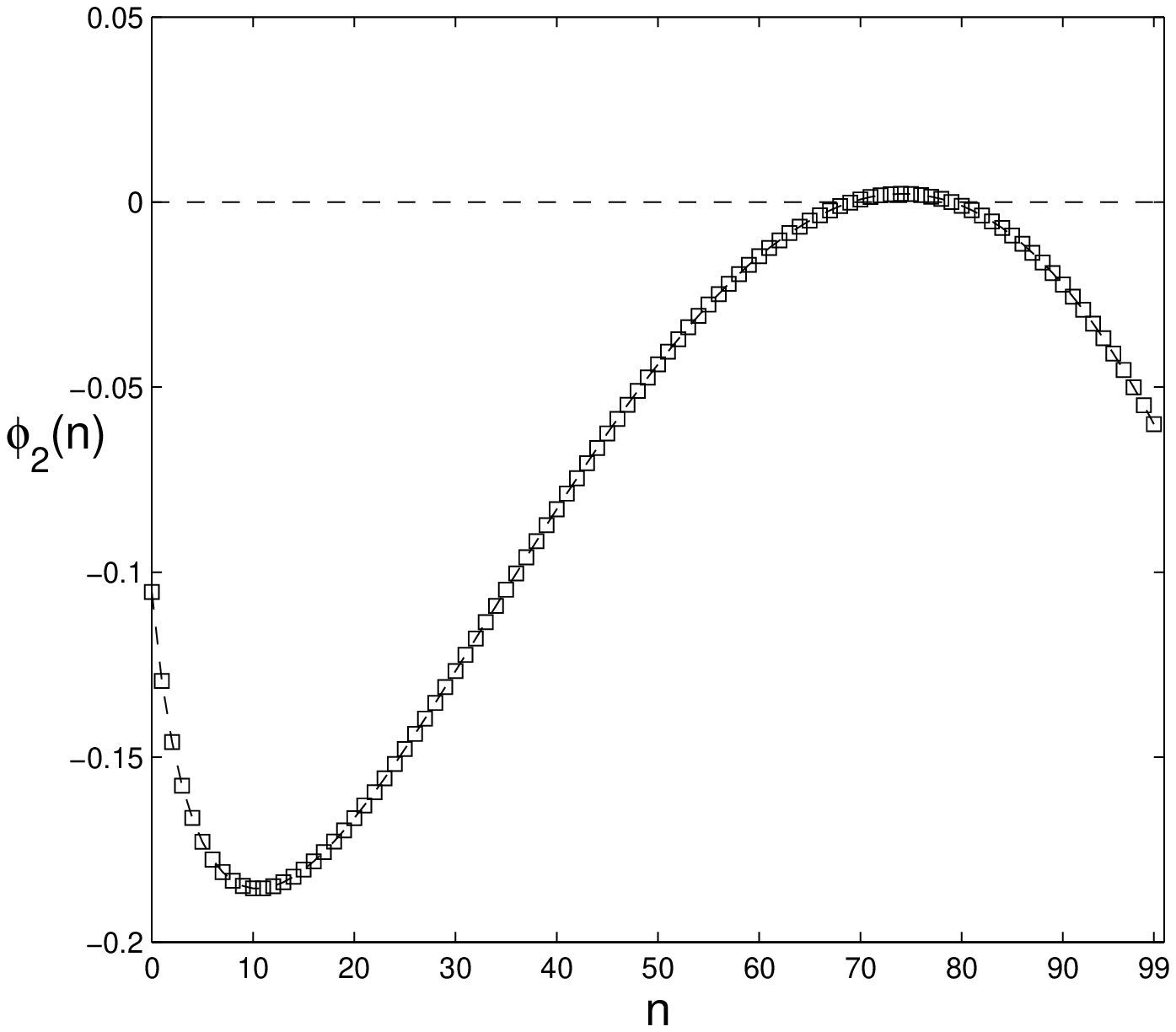}   
    \caption{$\phi_2(n)$ for $n\in[0,99]$ and $\rho=4$.}   
    \label{f2}   
  \end{minipage}%   
  \begin{minipage}[t]{0.5\linewidth}   
    \centering   
    \includegraphics[width=0.9\textwidth]{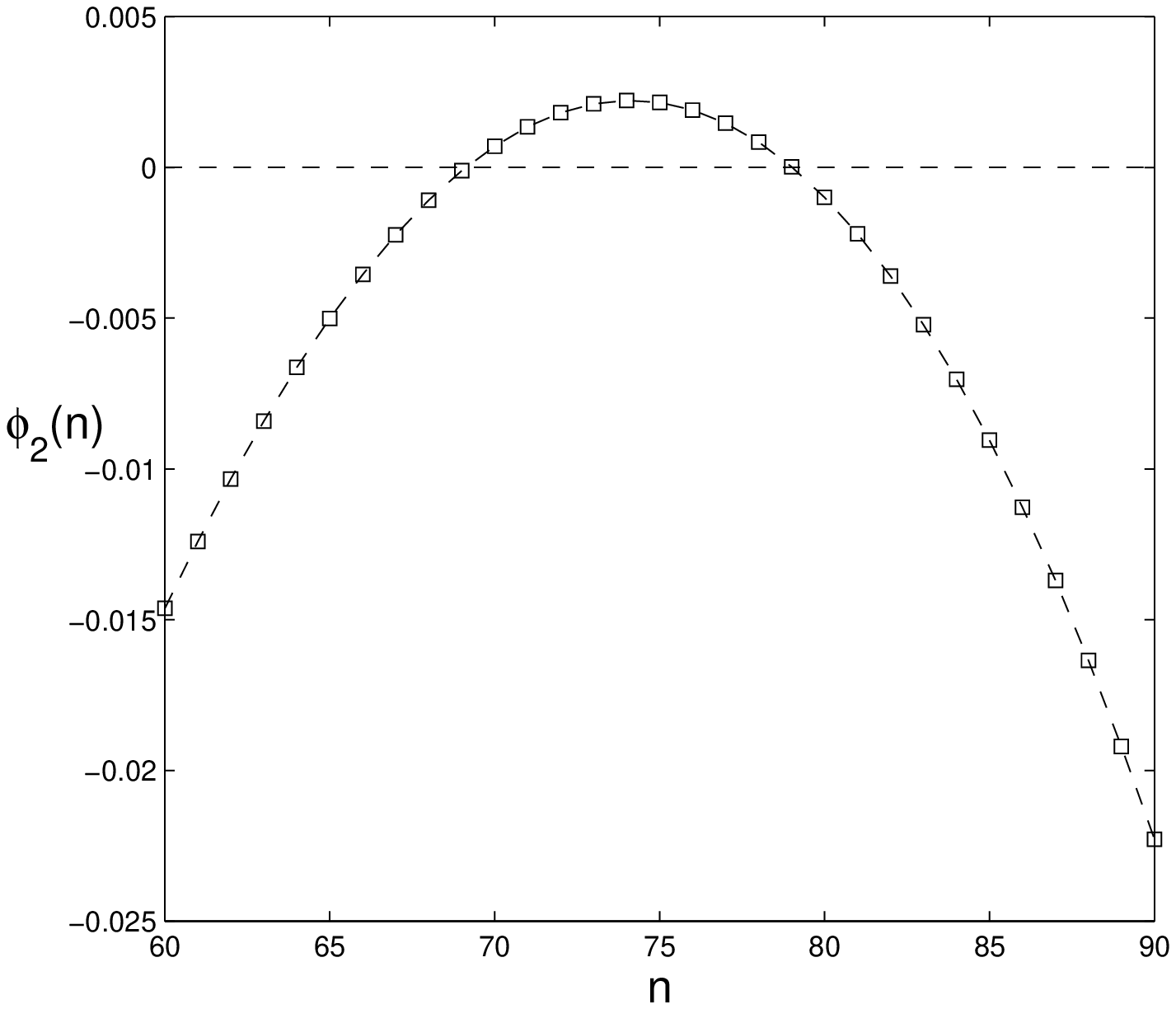}   
    \caption{$\phi_2(n)$ for $n\in[60,90]$ and $\rho=4$.}   
    \label{f3}   
  \end{minipage}   
\end{figure}

In Figure~\ref{f2} we retain $N=100$ and $\rho=4$, but now plot the second eigenvector $\phi_2(n)$. We see that typically $\phi_2(n)<0$, and its graph is approximately tangent to the $n$-axis near $n=74$. In Figure~\ref{f3} we ``blow up'' the region near $n=75$, plotting $\phi_2(n)$ for $n\in[60,90]$. Now we clearly see two sign changes, and these occur as $n$ changes from 69 to 70, and 79 to 80. Our asymptotic analysis suggests that near $X=3/4$ the eigenvector is proportional to $\mathrm{He}_2(\sqrt{\rho}X)=\rho X^2-1=4X^2-1$, which has zeros at $X=\pm 1/2$, and in view of (\ref{s31_X}) this corresponds to $n=75\pm 5$, which is again in excellent agreement with the exact results. In the range $n\in [60,90]$ Figure~\ref{f3} shows roughly a parabolic profile, as predicted by (\ref{s31_phijX}), but the larger picture in Figure~\ref{f2} again demonstrates that the $\xi$-scale result in (\ref{s31_phijk1}) must be used when $X$ is further away from $3/4$, as (\ref{s31_phijk1}) will, for example, predict the minimum value seen in Figure~\ref{f2}.

\begin{figure}   
  \begin{minipage}[t]{0.5\linewidth}  
    \centering   
    \includegraphics[width=0.9\textwidth]{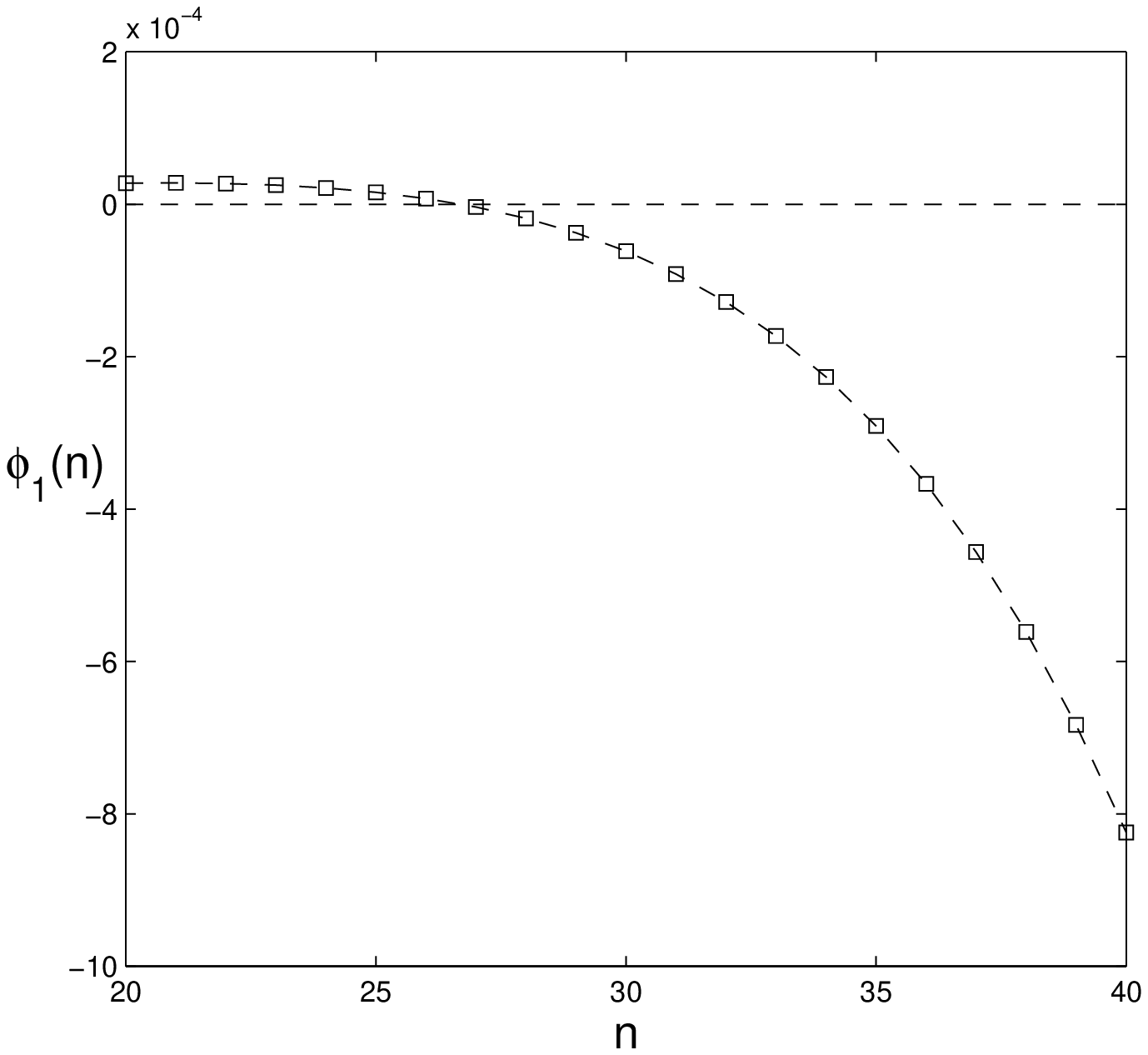}   
    \caption{$\phi_1(n)$ for $n\in[20,40]$ and $\rho=1$.}   
    \label{f5}   
  \end{minipage}%   
  \begin{minipage}[t]{0.5\linewidth}   
    \centering   
    \includegraphics[width=0.9\textwidth]{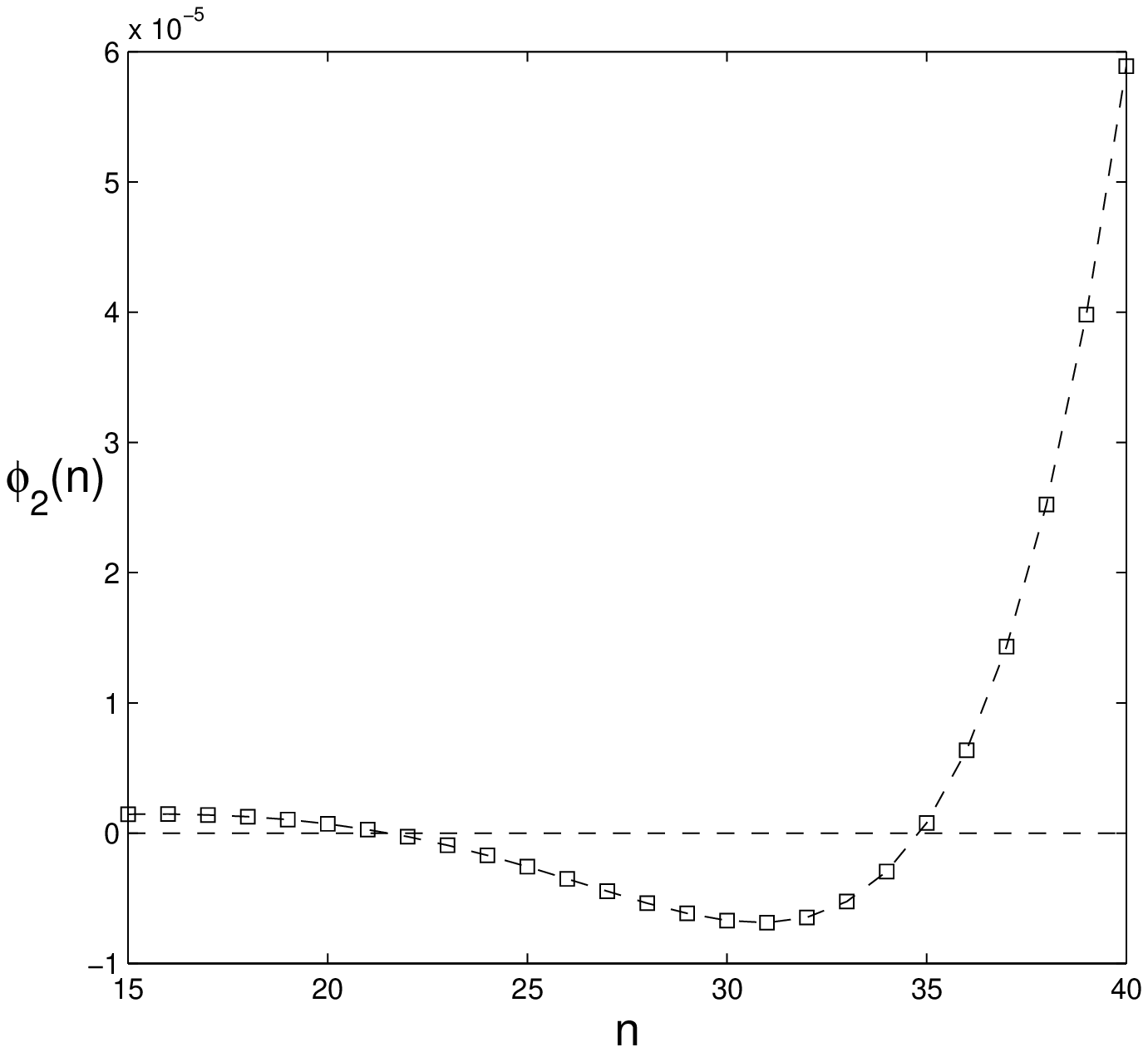}   
    \caption{$\phi_2(n)$ for $n\in[15,40]$ and $\rho=1$.}   
    \label{f6}   
  \end{minipage}   
\end{figure}

Next we consider $\rho=1$, maintaining $N=100$. Now the asymptotic result in the main range is in (\ref{s31_U}) and (\ref{s31_phieig}), and when $\rho=1$ we have $A=2^{1/3}$. In Figure~\ref{f5} we plot $\phi_1(n)$ ($j=1$) in the range $n\in[20,40]$, and we see a single sign change when $n$ increases from 26 to 27. For $j=1$ the asymptotic formula in (\ref{s31_phieig}) predicts a zero at $U=0$, which corresponds to $n=N^{2/3}A\approx 27.1$. In Figure~\ref{f6} we have $j=2$ and $n\in [15,40]$; there are two sign changes, between $n=21$ and 22, and $n=34$ and 35. Now $(1+4A^{-3})^{1/4}=3^{1/4}$ and the approximation in (\ref{s31_phieig}) has zeros where $\sqrt{3}\,U^2-1=0$, which corresponds to $n=N^{2/3}A\pm N^{1/2}3^{-1/4}$, leading to the numerical values of $n\approx 19.5$ and $n\approx 34.7$. 

\begin{figure}   
  \begin{minipage}[t]{0.5\linewidth}  
    \centering   
    \includegraphics[width=0.9\textwidth]{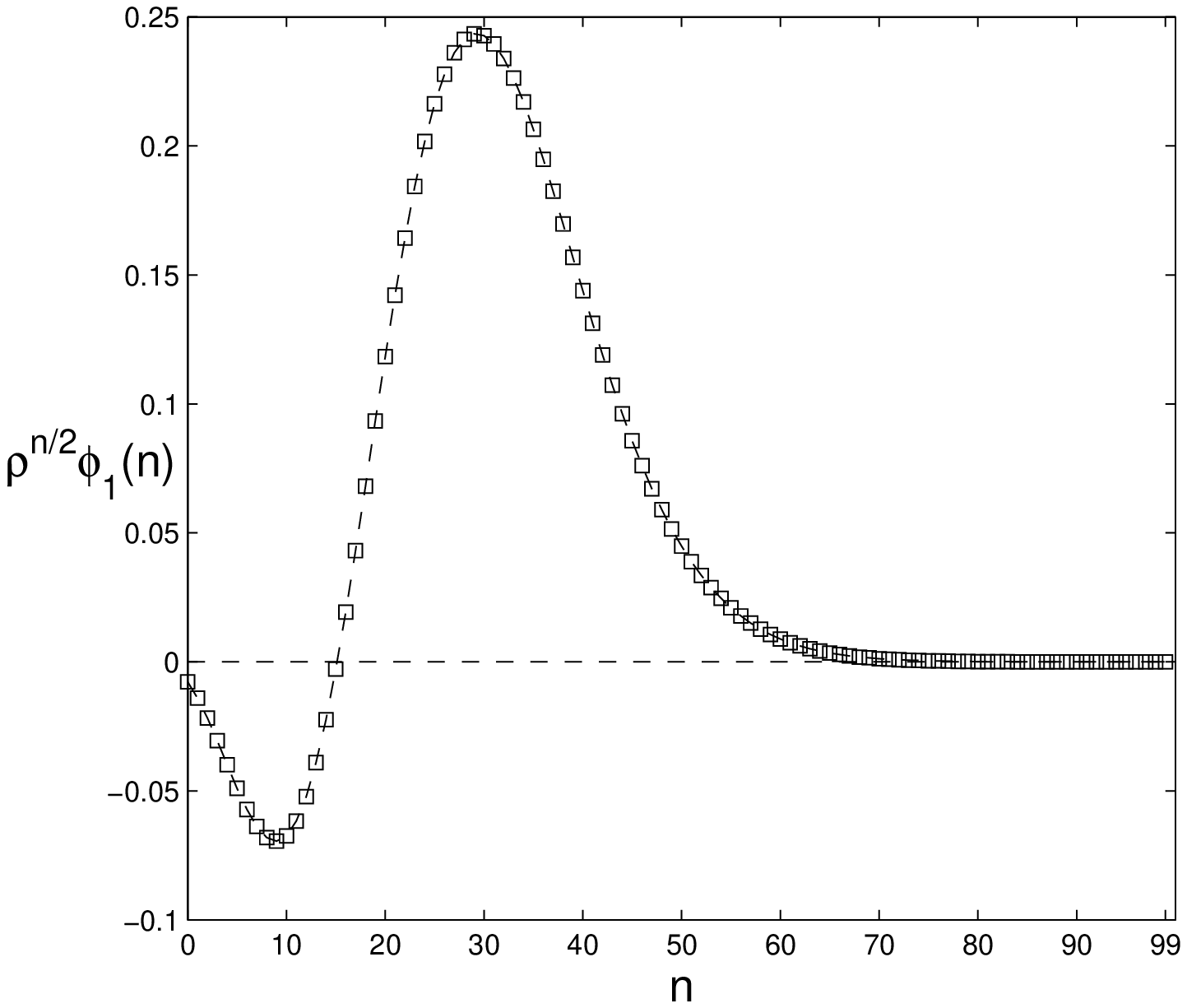}   
    \caption{$\rho^{n/2}\phi_1(n)$ for $n\in[0,99]$ and $\rho=0.25$.}   
    \label{f8}   
  \end{minipage}%   
  \begin{minipage}[t]{0.5\linewidth}   
    \centering   
    \includegraphics[width=0.9\textwidth]{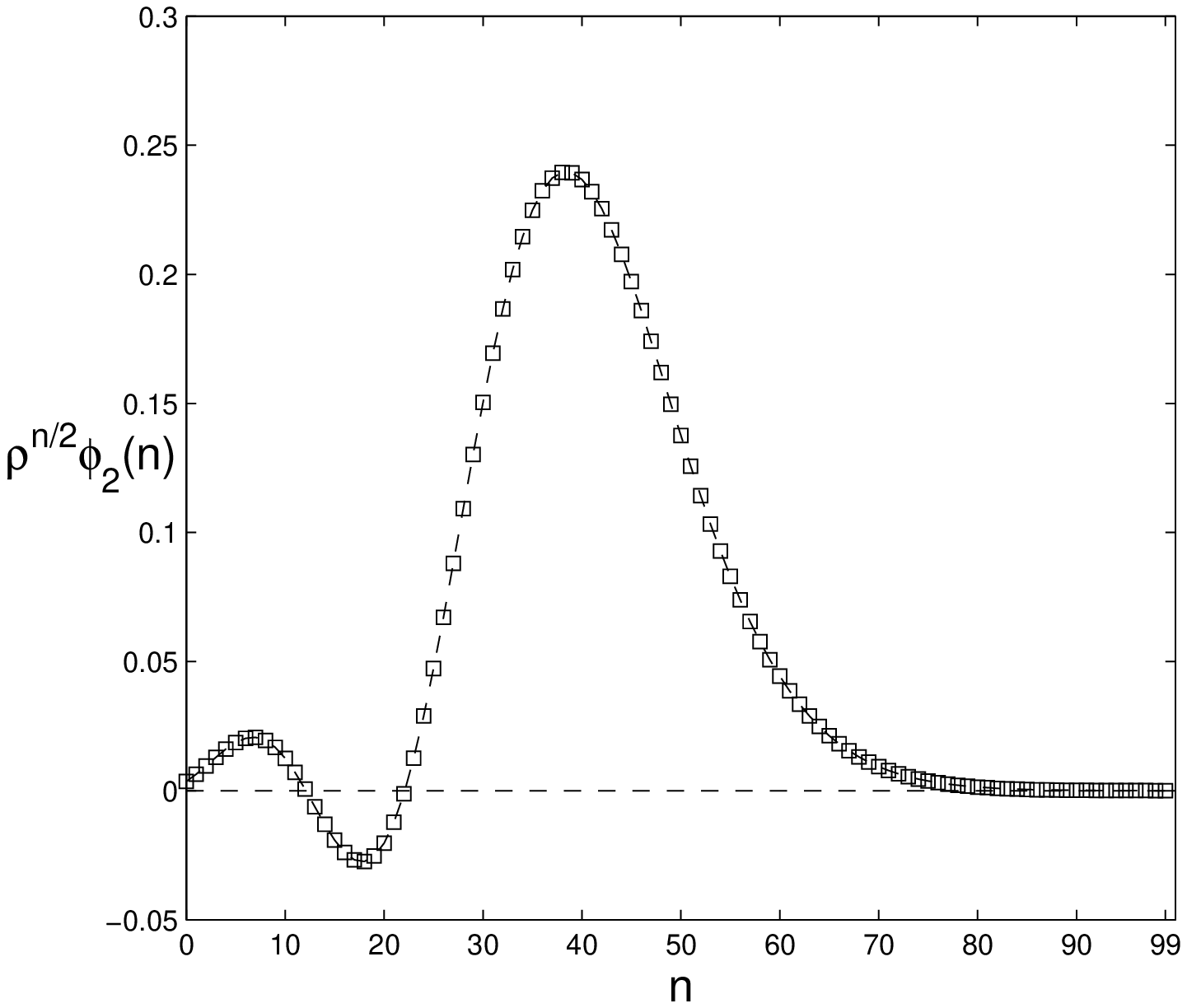}   
    \caption{$\rho^{n/2}\phi_2(n)$ for $n\in[0,99]$ and $\rho=0.25$.}   
    \label{f9}   
  \end{minipage}   
\end{figure}

Next we consider $\rho<1$. Now the main range is the $y$-scale result in (\ref{s31_phi})-(\ref{s31_ab}). For $\rho<1$ we now plot the ``symmetrized'' eigenvector(s) $\rho^{n/2}\phi_j(n)$. In Figures~\ref{f8} and \ref{f9} we always have $\rho=0.25$ and $N=100$. Figure~\ref{f8} has $\rho^{n/2}\phi_1(n)$ for $n\in [0,99]$ and shows a single sign change between $n=15$ and 16. The leading term in (\ref{s31_phi}) predicts a sign change at $y=0$, so that $n=\sqrt{N}/\sqrt{1-\sqrt{\rho}}=\sqrt{200}\approx 14.1$. Since $\Phi^{(1)}_1(0)=8\beta<0$ in (\ref{s31_Phi1}) (with (\ref{s31_ab})), the correction term would improve on the accuracy of the sign change prediction. Figure~\ref{f9} plots $\rho^{n/2}\phi_2(n)$ and shows two sign changes, between $n=12$ and 13, and $n=22$ and 23. Now the leading term in (\ref{s31_phi}) has zeros at $z=\pm 1$ which correspond to 
%$$n=\frac{\sqrt{N}}{\sqrt{1-\sqrt{\rho}}}\pm\frac{N^{3/8}}{\sqrt{2}(1-\sqrt{\rho})^{3/8}},$$
the numerical values $n\approx 8.99$ and $n\approx 19.30$ from (\ref{s31_ny}). We note that Figures~\ref{f8} and \ref{f9} show that the eigenvector becomes very small as $n$ increases toward 99, and this is indeed predicted by the expansion in (\ref{s31_phileading}), which applies on the $\xi=n/N$ scale.

We have thus shown that our asymptotic results predict quite well the qualitative properties of the eigenvectors $\phi_j(n)$ for $N$ large and moderate $j$. We get very good quantitative agreement when $\rho>1$, but less so when $\rho=1$ or $\rho<1$, for the moderately large value $N=100$. 

Next we consider the accuracy of our expansions for the eigenvalues $\nu_j$. In Table~\ref{table1} we take $\rho=0.25$ and increase $N$ from 10 to 100, and give the two-term, three-term, and four-term asymptotic approximations in (\ref{s31_nu}) for the zeroth eigenvalue $\nu_0$, along with the exact numerical result. We also give the relative errors, defined as $|\mathrm{num} -\mathrm{asy}|/\mathrm{num}$. We note that the leading term would simply give $(1-\sqrt{\rho})^2=0.25$, independent of $N$, which is not very accurate. Table~\ref{table1} shows that the correction terms in (\ref{s31_nu}) do lead to accurate approximations, and indeed when $N=100$ the four-term approximation is accurate to three significant figures. 

In Table~\ref{table2} we take $\rho=4>1$, in which case we only computed the leading term in (\ref{s31_nuj}), but now the eigenvalue index $j$ appears to leading order. Table~\ref{table2} gives the exact values of $\nu_0$ and $\nu_1$, along with the leading order approximations in (\ref{s31_nuj}), for $N$ increasing from 10 to 100. The agreement is again very good, with the relative errors decreasing to about $1\%$ when $N=100$.

We next do some numerical studies to see how rapidly the unconditional sojourn time density settles to its tail behavior, for problems with moderately large $N$. Here we compute $p(t)$ exactly (numerically), using (\ref{s2_fs_sum}) and (\ref{s2_fs_probN}) in (\ref{s1_pt}), and compare the result to the approximation in (\ref{s2_fs_ptsim}), which uses only the zeroth eigenvalue $\nu_0$. Tables~\ref{table3}-\ref{table5} compare the exact $-t^{-1}\log[p(t)]$ to the corresponding approximation from (\ref{s2_fs_ptsim}), and we note that both must approach $\nu_0$ as $t\to\infty$. Our asymptotic analysis predicts that for $N\to\infty$ and $\rho<1$, the $j=0$ term in (\ref{s2_fs_sum}) should dominate for times $t\gg O(N^{3/4})$, while if $\rho-1=O(N^{-1/3})$ or $\rho>1$, the $j=0$ term dominates for times $t\gg O(N)$. 

In Table~\ref{table3} we take $\rho=0.25$ and let $N=10$ and $20$. For $N=10$, the largest eigenvalue is $\nu_0\approx 0.5041$, which we list in the last row of the table, and the second largest eigenvalue is $\nu_1\approx 0.6290$. For $N=20$, we have $\nu_0\approx 0.4284$ and $\nu_1\approx 0.4982$. In Table~\ref{table4} we take $\rho=1$. Now $\nu_0\approx 0.2573$ and $\nu_1\approx 0.4642$ when $N=10$, and $\nu_0\approx 0.1635$ and $\nu_1\approx 0.2638$ when $N=20$. In Table~\ref{table5} we take $\rho=4$ and $\nu_0\approx 0.1312$, $\nu_1\approx 0.5805$ when $N=10$, and $\nu_0\approx 0.0662$, $\nu_1\approx 0.2784$ when $N=20$. 

Tables~\ref{table3}-\ref{table5} show that the approximation resulting from (\ref{s2_fs_ptsim}) is quite accurate, though it may take fairly large times before the exact and approximate results ultimately reach the limit $\nu_0$. The relative errors improve as we go from $\rho=0.25$ to $\rho=1$ to $\rho=4$, which is again consistent with our asymptotic analysis, as when $\rho<1$ there is the most coalescence (for $N\to\infty$) of the eigenvalues, making it hard to distinguish $\nu_0$ from the others.

\newpage
\vspace{2in}

\begin{center}
\includegraphics[width=0.45\textwidth]{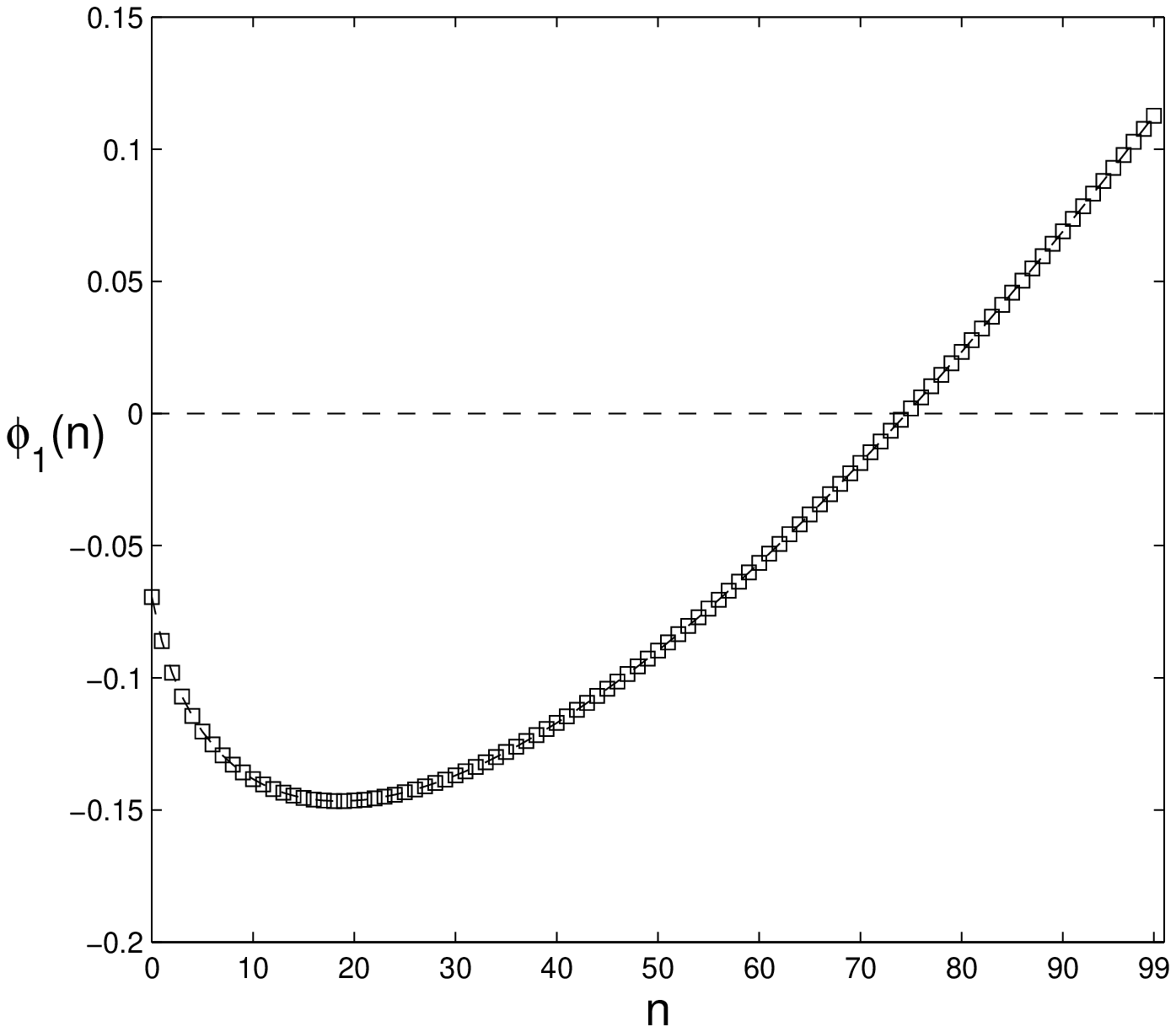}
\caption{$\phi_1(n)$ for $n\in[0,99]$ and $\rho=4$.} \label{f1}
\end{center}

\vspace{1.5in}

\begin{table}[H]
\centering
\caption{The eigenvalue $\nu_0$ with $\rho=0.25$.}\label{table1}
\vspace{4mm}
\begin{tabular} {|c||c||c|c||c|c||c|c|}       \hline
$N$		&	\tabincell{c}{$\nu_0$\\[-3pt] (Exact)}	& \tabincell{c}{2-term\\[-3pt] Approx.}		& \tabincell{c}{Relative\\[-3pt] Error}		& \tabincell{c}{3-term\\[-3pt] Approx.}		&	\tabincell{c}{Relative\\[-3pt] Error}		&	\tabincell{c}{4-term\\[-3pt] Approx.}		&	\tabincell{c}{Relative\\[-3pt] Error}  \\\hline
10&	0.5041&	0.4736&	6.0E-02&	0.5265&	4.4E-02&	0.4968&	1.5E-02\\\hline
%20&	0.4284&	0.4081&	4.74E-02&	0.4395&	2.60E-02&	0.4247&	8.69E-03\\\hline
30&	0.3948&	0.3791&	4.0E-02&	0.4023&	1.9E-02&	0.3924&	6.2E-03\\\hline
%40&	0.3749&	0.3618&	3.48E-02&	0.3805&	1.50E-02&	0.3731&	4.78E-03\\\hline
50&	0.3613&	0.3500&	3.1E-02&	0.3658&	1.3E-02&	0.3599&	3.9E-03\\\hline
%60&	0.3513&	0.3413&	2.84E-02&	0.3551&	1.08E-02&	0.3501&	3.26E-03\\\hline
70&	0.3435&	0.3345&	2.6E-02&	0.3468&	9.5E-03&	0.3426&	2.8E-03\\\hline
%80&	0.3373&	0.3291&	2.44E-02&	0.3402&	8.55E-03&	0.3365&	2.46E-03\\\hline
%90&	0.3321&	0.3245&	2.29E-02&	0.3347&	7.75E-03&	0.3314&	2.18E-03\\\hline
100&	0.3278&	0.3207&	2.2E-02&	0.3301&	7.1E-03&	0.3271&	2.0E-03\\\hline
\end{tabular}\\
\end{table}

\begin{table}[H]
\centering
\caption{The eigenvalues $\nu_0$ and $\nu_1$ with $\rho=4$.}\label{table2}
\vspace{4mm}
\begin{tabular} {|c||c|c|c||c|c|c|}       \hline
$N$		&	\tabincell{c}{$\nu_0$\\[-3pt] (Exact)}		& \tabincell{c}{$\nu_0$\\[-3pt] (Approx.)}		& \tabincell{c}{Relative\\[-3pt] Error}		&	\tabincell{c}{$\nu_1$\\[-3pt] (Exact)}		& \tabincell{c}{$\nu_1$\\[-3pt] (Approx.)}		&	\tabincell{c}{Relative\\[-3pt] Error} \\\hline
10&0.1312&0.1333&1.7E-02&0.5805&0.5333&8.1E-02\\\hline
%20&0.0662&0.0667&7.77E-03&0.2784&0.2667&4.22E-02\\\hline
30&0.0442&0.0444&5.1E-03&0.1830&0.1778&2.8E-02\\\hline
%40&0.0332&0.0333&3.79E-03&0.1362&0.1333&2.13E-02\\\hline
50&0.0266&0.0267&3.0E-03&0.1085&0.1067&1.7E-02\\\hline
%60&0.0222&0.0222&2.51E-03&0.0902&0.0889&1.42E-02\\\hline
70&0.0190&0.0190&2.1E-03&0.0771&0.0762&1.2E-02\\\hline
%80&0.0166&0.0167&1.87E-03&0.0674&0.0667&1.06E-02\\\hline
%90&0.0148&0.0148&1.66E-03&0.0598&0.0593&9.46E-03\\\hline
100&0.0133&0.0133&1.5E-03&0.0538&0.0533&8.5E-03\\\hline
\end{tabular}\\
\end{table}

\begin{table}[H]
\centering
\caption{The tail approximation of $p(t)$ with $\rho=0.25$.}\label{table3}
\vspace{4mm}
\begin{tabular} {|c||c|c|c||c|c|c|}       \hline
&\multicolumn{3}{|c||}{$\rho=0.25,\;N=10$}  &\multicolumn{3}{c|}{$\rho=0.25,\;N=20$}  \\ \hline
$t$	& \tabincell{c}{$-\log[p(t)]/t$\\[-3pt] (Exact)} 	& \tabincell{c}{$-\log[p(t)]/t$\\[-3pt] (Approx.)}	& \tabincell{c}{Relative\\[-3pt] Error}	& \tabincell{c}{$-\log[p(t)]/t$\\[-3pt] (Exact)} 	& \tabincell{c}{$-\log[p(t)]/t$\\[-3pt] (Approx.)}		&\tabincell{c}{Relative\\[-3pt] Error} \\\hline
5	&0.8186	&1.0204 &24.66\%	&0.8078	&1.1992	&48.462\% \\\hline
10	&0.7141	&0.7623	&6.74\%		&0.6947	&0.8138	&17.16\%	\\\hline
15	&0.6762	&0.6599	&2.47\%		&0.6854	&0.6339	&8.12\%  \\\hline
20	&0.6268	&0.6332	&1.02\%		&0.5953	&0.6211	&4.33\%	\\\hline
30	&0.5890	&0.5902	&0.20\%		&0.5488	&0.5569	&1.47\%	\\\hline
%40	&0.5684	&0.5687	&0.05\%		&0.5219	&0.5248	&0.56\%	\\\hline
50	&0.5557	&0.5558	&0.01\%		&0.5044	&0.5055	&0.23\%	\\\hline
100	&0.5299	&0.5299	&1.07E-07	&0.4670	&0.4670	&3.69E-05 \\\hline	
%500	&0.5093	&0.5093	& $<$1E-12	&0.4361	&0.4361	& $<$1E-12	\\\hline
$\infty$	&0.5041	&0.5041		&--	&0.4284	&0.4284		&--\\\hline
\end{tabular}
\end{table}

\begin{table}[H]
\centering
\caption{The tail approximation of $p(t)$ with $\rho=1$.}\label{table4}
\vspace{4mm}
\begin{tabular} {|c||c|c|c||c|c|c|}       \hline
&\multicolumn{3}{|c||}{$\rho=1,\;N=10$}  &\multicolumn{3}{c|}{$\rho=1,\;N=20$}  \\ \hline
$t$	& \tabincell{c}{$-\log[p(t)]/t$\\[-3pt] (Exact)} 	& \tabincell{c}{$-\log[p(t)]/t$\\[-3pt] (Approx.)}	& \tabincell{c}{Relative\\[-3pt] Error}	& \tabincell{c}{$-\log[p(t)]/t$\\[-3pt] (Exact)} 	& \tabincell{c}{$-\log[p(t)]/t$\\[-3pt] (Approx.)}		&\tabincell{c}{Relative\\[-3pt] Error} \\\hline
5	&0.6030	&0.6664	&10.52\% &0.5800	&0.7295	&25.77\%	\\\hline
10	&0.4515	&0.4618	&2.30\%	&0.4070	&0.4465	&9.71\%	\\\hline
15	&0.3913	&0.3937	&0.61\%	&0.3371	&0.3521	&4.45\%	\\\hline
20	&0.3589	&0.3596	&0.18\%	&0.2984	&0.3050	&2.21\%	\\\hline
30	&0.3254	&0.3255	&0.02\%	&0.2562	&0.2578	&0.61\%	\\\hline
50	&0.2982	&0.2982	&1.72E-06 &0.2435	&0.2443	&0.33\%	\\\hline
100	&0.2778	&0.2778	&2.96E-11 &0.2338	&0.2342	&0.18\%	\\\hline
$\infty$	&0.2573	&0.2573	&--	&0.1635	&0.1635	&--	\\\hline
\end{tabular}
\end{table}

\begin{table}[H]
\centering
\caption{The tail approximation of $p(t)$ with $\rho=4$.}\label{table5}
\vspace{4mm}
\begin{tabular} {|c||c|c|c||c|c|c|}       \hline
&\multicolumn{3}{|c||}{$\rho=4,\;N=10$}  &\multicolumn{3}{c|}{$\rho=4,\;N=20$}  \\ \hline
$t$	& \tabincell{c}{$-\log[p(t)]/t$\\[-3pt] (Exact)} 	& \tabincell{c}{$-\log[p(t)]/t$\\[-3pt] (Approx.)}	& \tabincell{c}{Relative\\[-3pt] Error}	& \tabincell{c}{$-\log[p(t)]/t$\\[-3pt] (Exact)} 	& \tabincell{c}{$-\log[p(t)]/t$\\[-3pt] (Approx.)}		&\tabincell{c}{Relative\\[-3pt] Error} \\\hline
5	&0.5398	&0.5415	&0.32\%	&0.3383	&0.3387	&0.14\%	\\\hline
10	&0.3362	&0.3363	&0.03\%	&0.2024	&0.2024	&0.01\%	\\\hline
15	&0.2679	&0.2680	&2.4E-05	&0.1570	&0.1570	&1.5E-05	\\\hline
20	&0.2338	&0.2338	&2.2E-06	&0.1343	&0.1343	&1.6E-06	\\\hline
30	&0.1996	&0.1996	&1.9E-08	&0.1207	&0.1207	&1.7E-07	\\\hline
50	&0.1722	&0.1722	&1.6E-12	&0.1002	&0.1002	&2.1E-10	\\\hline
100	&0.1517	&0.1517	&$<$1.0E-12	&0.0934	&0.0934	&2.6E-12	\\\hline
$\infty$	&0.1312	&0.1312	&--	&0.0662	&0.0662	&--	\\\hline
\end{tabular}
\end{table}


\newpage

\begin{thebibliography}{2}
\bibitem{KLE_A} L. Kleinrock, {\it Analysis of a time-shared processor}, Naval Research Logistics Quarterly 11 (1964), 59-73.\\[-20pt]
\bibitem{KLE_T} L. Kleinrock, {\it Time-shared systems: A theoretical treatment}, J. ACM 14 (1967), 242-261.\\[-20pt]
\bibitem{HEY} D. P. Heyman, T. V. Lakshman, and A. L. Neidhardt, {\it A new method for analysing feedback-based protocols with applications to engineering Web traffic over the Internet}, Proc. ACM Sigmetrics (1997), 24-38.\\[-20pt]
\bibitem{MAS} L. Massouli\'e and J. W. Roberts, {\it Bandwidth sharing: Objectives and algorithms}, Proc. IEEE Infocom, New York, NY, USA (1999), 1395-1403.\\[-20pt]
\bibitem{NAB} M. Nabe, M. Murata, and H. Miyahara, {\it Analysis and modeling of World Wide Web traffic for capacity dimensioning of Internet access lines}, Perf. Evaluation, 34 (1998), 249-271.\\[-20pt]
\bibitem{MIT} D. Mitra and J. A. Morrison, {\it Asymptotic expansions of moments of the waiting time in closed and open processor-sharing systems with multiple job classes}, Adv. in Appl. Probab. 15 (1983), 813-839.\\[-20pt]
\bibitem{MOR_H} J. A. Morrison and D. Mitra, {\it Heavy-usage asymptotic expansions for the waiting time in closed processor-sharing systems with multiple classes}, Adv. in Appl. Probab. 17 (1985), 163-185.\\[-20pt]
\bibitem{MOR_A} J. A. Morrison, {\it Asymptotic analysis of the waiting-time distribution for a large closed processor-sharing system}, SIAM J. Appl. Math. 46 (1986), 140-170.\\[-20pt]
\bibitem{MOR_M} J. A. Morrison, {\it Moments of the conditioned waiting time in a large closed processor-sharing system}, Stochastic Models 2 (1986), 293-321.\\[-20pt]
\bibitem{MOR_C87} J. A. Morrison, {\it Conditioned response-time distribution for a large closed processor-sharing system in very heavy usage}, SIAM J. Appl. Math. 47 (1987), 1117-1129.\\[-20pt]
\bibitem{MOR_C88} J. A. Morrison, {\it Conditioned response-time distribution for a large closed processor-sharing system with multiple classes in very heavy usage}, SIAM J. Appl. Math. 48 (1988), 1493-1509.\\[-20pt]
\bibitem{SEV} K. C. Sevcik and I. Mitrani, {\it The distribution of queueing network states at input and output instants}, J. ACM 28 (1981), 358-371.\\[-20pt]
\bibitem{POL} F. Pollaczek, {\it La loi d'attente des appels t\'el\'ephoniques}, C. R. Acad. Sci. Paris 222 (1946), 353-355.\\[-20pt]
\bibitem{COH} J. W. Cohen, {\it On processor sharing and random service (Letter to the editor)}, J. Appl. Prob. 21 (1984), 937.\\[-20pt]
\bibitem{mag} W. Magnus, F. Oberhettinger, and R. P. Soni, Formulas and Theorems for the Special Functions of Mathematical Physics, Springer-Verlag, New York, 1966.

\end{thebibliography}
\end{document}